\RequirePackage{ifpdf}
\ifpdf 
\documentclass[pdftex]{sigma}
\else
\documentclass{sigma}
\fi

\let\Bbb=\mathbb

\def\XX{{
{{\mbox{\tiny$\times$}}\atop{
\mbox{\tiny$\times$}}}}}
\newcommand{\ORD}[1]{\XX{#1}\XX}

\newcommand{\tr}{\,{\rm tr}\,}

\def\e{e}

\def\RR{{\Bbb R}}

\def\mod{\mathop{\mbox{mod}}\nolimits}

\begin{document}

\numberwithin{equation}{section}

\allowdisplaybreaks

\renewcommand{\PaperNumber}{066}

\FirstPageHeading

\renewcommand{\thefootnote}{$\star$}

\ShortArticleName{Teichm\"uller Theory of Bordered Surfaces}

\ArticleName{Teichm\"uller Theory of Bordered Surfaces\footnote{This paper is a
contribution to the Vadim Kuznetsov Memorial Issue `Integrable
Systems and Related Topics'. The full collection is available at
\href{http://www.emis.de/journals/SIGMA/kuznetsov.html}{http://www.emis.de/journals/SIGMA/kuznetsov.html}}}

\Author{Leonid O. CHEKHOV}

\AuthorNameForHeading{L.O. Chekhov}

\Address{Steklov Mathematical Institute, Moscow, Russia}
\Email{\href{mailto:chekhov@mi.ras.ru}{chekhov@mi.ras.ru}}

\Address{Institute for Theoretical and Experimental Physics, Moscow, Russia}

\Address{Poncelet Laboratoire International Franco--Russe, Moscow, Russia}

\Address{Concordia University, Montr\'eal, Quebec, Canada}

\ArticleDates{Received January 05, 2007, in f\/inal form April
28, 2007; Published online May 15, 2007}

\Abstract{We propose the graph description of Teichm\"uller theory of surfaces
with marked points on boundary components (bordered surfaces).
Introducing new parameters, we formu\-late this theory in terms of
hyperbolic geometry. We can then describe both classical and quantum
theories having the proper number of Thurston variables
(foliation-shear coordinates), mapping-class group invariance (both
classical and quantum), Poisson and quantum algebra of geodesic
functions, and classical and quantum braid-group relations. These
new algebras can be def\/ined on the double of the corresponding graph
related (in a novel way) to a double of the Riemann surface (which
is a Riemann surface with holes, not a smooth Riemann surface). We
enlarge the mapping class group allowing transformations relating
dif\/ferent Teichm\"uller spaces of bordered surfaces of the same
genus, same number of boun\-dary components, and same total number of
marked points but with arbitrary distributions of marked points
among the boundary components. We describe the classical and quantum
algebras and braid group relations for particular sets of geodesic
functions correspon\-ding to~$A_n$ and~$D_n$ algebras and discuss
brief\/ly the relation to the Thurston theory.}


\Classification{37D40; 53C22}

\renewcommand{\thefootnote}{\arabic{footnote}}
\setcounter{footnote}{0}

\section{Introduction}

Recent advances in the quantitative description of the Teichm\"uller spaces of hyperbolic struc\-tures were mainly based on
the graph (combinatorial) description of the corresponding spa\-ces~\mbox{\cite{Penn1,Fock1}}. The corresponding structures not only
provided a convenient coordinatization together with the mapping class group action, they proved to be especially useful
when describing sets of geodesic functions and the related Poisson and quantum structures \cite{ChF}. Combined with Thurston's
theory of measured foliations \cite{ThSh,HarPen}, it led eventually to the formulation of the quantum Thurston
theory~\cite{ChP}. The whole consideration was concerning Riemann surfaces with holes. A natural generalization of this
pattern consists in adding marked points on the boundary components. First, Kaufmann and Penner~\cite{KP} showed that the related
Thurston theory of measured foliations provides a nice combinatorial description of open/closed string diagrammatic. Second, if
approaching these systems from the algebraic viewpoint, one can associate a~cluster algebra  (originated in \cite{FZ} and applied
to bordered surfaces in \cite{FST}) to such a geometrical pattern.

\looseness=1
The aim of this paper is to provide a shear-coordinate description of Teichm\"uller spaces of bordered Riemann surfaces,
to construct the corresponding geodesic functions (cluster variables), and to investigate the Poisson and quantum relations satisf\/ied
by these functions in classical case or by the correspondent Hermitian operators in the quantum case.

\looseness=1
In Section~\ref{s:hyp}, we give a (presumably new) description of the Teichm\"uller space of bordered (or windowed) surfaces in the
hyperbolic geometry pattern using the graph technique supporting it by considering a simplest example of annulus with one marked point.
It turns out that adding each new window (a new marked point) increases the number of parameters by two resulting in adding a new
inversion relation to the set of the Fuchsian group generators. We explicitly formulate rules by which we can construct geodesic functions
(corresponding to components of a multicurve) using these coordinates; the only restriction we impose and keep throughout the paper is
the evenness condition: an even number of multicurve lines must terminate at each window.

\looseness=1
In Section~\ref{s:algebra}, we construct algebras of geodesic functions postulating the Poisson relations on the level of the (old and new)
shear coordinates of the Teichm\"uller space. We construct f\/lip morphisms and the corresponding mapping class group transformations and
f\/ind that in the bordered surfaces case we can enlarge this group allowing transformations that permute marked points on one of the
boundary components or transfer marked points from one component to another thus establishing isomorphisms between all the Teichm\"uller
spaces of surfaces of the same genus, same number of boundary components, and the same total number of marked points. In the same section,
we describe geodesic algebras corresponding (in the cluster terminology, see~\cite{FST}) to~$A_n$ and~$D_n$ systems. Whereas the
$A_n$-algebras have been known previously as algebras of geodesics on Riemann surfaces of higher genus \cite{NR,NRZ} (their graph
description in the case of higher-genus surfaces with one or two holes see in \cite{ChF2}) or as the algebra of Stockes parameters
\cite{DM,Ugaglia}, or as the algebra of upper-triangular matrices \cite{Bondal}, the $D_n$-algebras seem to be of a new sort. Using the
new type of the mapping class group transformations, we prove the braid group relations for all these algebras.

Section~\ref{ss:q} is devoted to quantization. We begin with a brief accounting of the quantization procedure from \cite{ChF} coming then to
the quantum geodesic operators and to the corresponding quantum algebras. Here, again, the quantum $D_n$-algebras seem to be of
a new sort, and we prove the Jacobi identities for them in the abstract setting without appealing to geometry. We also construct the
quantum braid group action in this section.

In Section~\ref{s:double}, we describe multicurves and related foliations for bordered surfaces, that is, we construct elements of Thurston's theory.
There we also explicitly construct the relevant doubled Riemann surface, which, contrary to what one could expect, is itself a Riemann surface
with holes (but without windows). We transfer the notion of multicurves to this doubled surface. Note, however, that the
new mapping class group transformations, while preserving the multicurve structure on the original bordered surface, change the
topological type of the doubled Riemann surface, which can therefore be treated only as an auxiliary, not basic, element of the
construction. Using this double, we can nevertheless formulate the basic statement similar to that in \cite{ChP}, that is,
that in order to obtain a self-consistent theory that is continuous at Thurston's boundary, we must set into the
correspondence to a multicurve the sum of lengths of its constituting geodesics (the sum of proper length operators in
the quantum case). In the same Section~\ref{s:double}, we describe elements of Thurston's theory of measured foliations for
bordered Riemann surfaces and the foliation-shear coordinate changings under the ``old'' and ``new'' mapping class group
transformations.

We discuss some perspectives of the proposed theory in the concluding section.

We tried to make the presentation as explicit and attainable as possible, so there are many f\/igures in the text.

\newpage

\section{Graph description and hyperbolic geometry\label{s:hyp}}

\subsection{Hyperbolic geometry and inversion relation\label{ss:graph}}

\subsubsection{Graph description for nonbordered Riemann surfaces}

Recall the graph description of the Fuchsian group, or the
fundamental group of the surface~$\Sigma_{g,s}$, which is a discrete
f\/initely generated subgroup of $PSL(2,\RR)$. We consider a spine
$\Gamma_{g,s}$ corresponding to the Riemann surface $\Sigma_{g,s}$
with $g$ handles and $s$ boundary components (holes). The spine, or
fatgraph $\Gamma_{g,s}$ is a connected graph that can be drawn
without self-intersections on $\Sigma_{g,s}$, has all vertices of
valence three, has a prescribed cyclic ordering of labeled edges
entering each vertex, and is a maximum graph in the sense that after
cutting along all its edges, the Riemann surface decomposes into the
set of polygons (faces) such that each polygon contains exactly one
hole (and becomes simply connected after plumbing this hole). Since
a graph must have at least one face, we can therefore describe only
Riemann surfaces with at least one hole, $s>0$. The hyperbolicity
condition also implies $2g-2+s>0$. We do not impose restrictions,
for instance, we allow edges to start and terminate at the same
vertex, allow two vertices to be connected with more than one edge,
etc. We however demand a spine to be a cell complex, that is, we do
not allow loops without vertices.

Then, we can establish a 1-1 correspondence between elements of the
Fuchsian group and closed paths in the spine starting and
terminating at the same directed edge. Since the terms in the matrix
product depend on the turns in vertices (see below), it is not
enough to f\/ix just a starting vertex. To construct an element of the
Fuchsian group $\Delta_{g,s}$, we select a directed edge (one and
the same for all the elements; see the example in
Fig.~\ref{fi:annulus} where it is indicated by a short fat arrow),
then move along edges and turns of the graph without backtracking
and eventually turn back to the selected directed edge\footnote{If
the last edge was not the selected one but its neighboring edge, the
very last move is turning to the selected edge, that is, we add
either $R$- or $L$-matrix; if the last edge coincides with the
selected edge, we do not make the last turn through the angle
$2\pi$. This results in the ambiguity by multiplication by the
matrix $-\hbox{Id}=R^3$, but it is inessential as we deal with the
projective transformations.}.

We associate with the $\alpha$th edge of the graph the real $Z_\alpha$ and set~\cite{Fock1}
the matrix of the M\"obius transformation
\begin{gather}
\label{XZ}
X_{Z_\alpha}=\left(
\begin{array}{cc} 0 & -\e^{Z_\alpha/2}\\
                \e^{-Z_\alpha/2} & 0\end{array}\right)
\end{gather}
each time the path homeomorphic to a geodesic $\gamma$ passes through the $\alpha$th edge.

We also introduce the ``right'' and ``left'' turn matrices
to be set in the proper place when a~path makes the corresponding turn,
\begin{gather*}
R=\left(\begin{array}{cc} 1 & 1\\ -1 & 0\end{array}\right), \qquad
L= R^2=\left(\begin{array}{cc} 0 & 1\\ -1 &
-1\end{array}\right),
\end{gather*}
and def\/ine the related operators $R_Z$ and $L_Z$,
\begin{gather*}
R_Z\equiv RX_Z=\left(\begin{array}{cc}
                \e^{-Z/2}&-\e^{Z/2}\\
                     0   &\e^{Z/2}
                     \end{array}\right),\\
L_Z\equiv LX_Z=\left(\begin{array}{cc}
                \e^{-Z/2}&   0\\
                 -\e^{-Z/2}&\e^{Z/2}
                     \end{array}\right).
\end{gather*}

An element of a Fuchsian group has then the structure
\[
P_{\gamma}=LX_{Z_n}RX_{Z_{n-1}}\cdots RX_{Z_2}RX_{Z_1},
\]
and the corresponding {\em geodesic function}
\begin{gather}
\label{G}
G_{\gamma}\equiv \tr P_\gamma=2\cosh(\ell_\gamma/2)
\end{gather}
is related to the actual length $\ell_\gamma$ of the closed
geodesic on the Riemann surface.

\subsubsection{Generalization to the bordered surfaces case}

We now introduce a new object, the {\em marking} pertaining to boundary components.
Namely, we assume that we have not just boundary components but allow some of them to
carry a f\/inite number (possibly zero) of marked points. We let $\delta_i$, $i=1,\dots,s$,
denote the corresponding number of marked points for the $i$th boundary component. Geometrically,
we assume these points to lie on the absolute, that is, instead of associating a closed geodesic
to the boundary component in nonmarked case, we associate to an $i$th boundary component
a collection comprising~$\delta_i$ inf\/inite geodesic curves connecting neighbor (in the sense of
the surface orientation) marked points on the absolute (can be the same point if $\delta_i=1$)
in the case where $\delta_i>0$. All these additional geodesic curves are disjoint with each other
and disjoint with any closed geodesic on the Riemann surface.
In~\cite{KP}, these curves were called {\em windows}. We denote
the corresponding windowed surface $\Sigma_{g,\delta}$, where
\begin{gather}
\label{delta}
\delta=\{\delta_1,\dots,\delta_s\}
\end{gather}
is the multiindex counting marked points on the boundary components ($\delta_i$ can be zero)
whereas $s$ is the number of
boundary components. We call such Riemann surfaces the windowed, or bordered Riemann surfaces.

Restrictions on $g$, $s$, and the number of marked points $\#\delta$ can be uniformly written
as $s>0$ and $2g-2+s+\left[\frac{\#\delta+1}{2}\right]>0$, that is, we allow two new cases
$g=0$, $s=1$, $\#\delta\ge 3$ and $g=0$, $s=2$, $\#\delta\ge 1$.

We want now to generalize the graph setting to the case where we have boundary components with
marked points. However, as the example below shows, in order to def\/ine inambiguously the
corresponding hyperbolic geometry, when introducing a marked point on the boundary, we must
simultaneously introduce {\em one more} additional parameter. This is because, as we shall demonstrate,
introducing a marked point adds a new inversion relation that preserves the orientation
but not the surface itself,
that is, we
invert a part of the Riemann surface through a boundary curve without taking care on
what happen to the (remaining) part of the surface because,
in our description, this (outer) part is irrelevant. Such an
inversion relation leaves invariant the new added geodesic, that is, the
window. However, there is
a {\em one-parameter family} of such inversions for every window, and in order to f\/ix the ambiguity we must indicate explicitly
which point on the new geodesic is stable w.r.t.\ such an inversion. Recall that because of
orientation preservation, two ends (on the absolute) of this new geodesic must be interchanged by
the inversion relation; it is therefore a unique point that is stable.

\begin{figure}[tb]
\centerline{\includegraphics{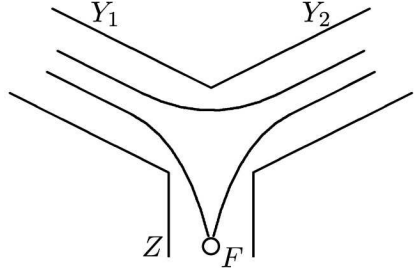}}
\caption{Part of the
spine with the pending vertex. The
variable $Z$ corresponds to the respective pending edge. Two
types of geodesic lines are shown in the f\/igure: one that does not
come to the edge $Z$ is parameterized in the standard way, the
other undergoes the inversion with the matrix $F$~(\ref{F}).}
\label{fi:corner}
\end{figure}

We describe this situation by considering {\em new types of graphs} with pending vertices.
Assume that we have a part of graph having the structure as in Fig.~\ref{fi:corner}.

Then, if a geodesic line comes to a pending vertex, it undergoes the
{\em inversion},
which stems to that we insert the {\em inversion matrix} $F$,
\begin{gather}
\label{F}
F=\left(\begin{array}{cc} 0 & 1\\ -1 & 0\end{array}\right),
\end{gather}
into the corresponding string of $2\times2$-matrices. For example, a part of geodesic function
in Fig.~\ref{fi:corner} that is inverted reads
\[
\cdots X_{Y_1}LX_Z FX_Z LX_{Y_2}\cdots,
\]
whereas the other geodesic that does not go to the pending vertex reads merely
\[
\cdots X_{Y_1} RX_{Y_2}\cdots.
\]
We call this new relation the inversion relation, and the
inversion element is itself an element of $PSL(2,\RR)$. We also call the edge terminating at a pending
vertex the {\em pending} edge.

Note the simple relation\footnote{In particular, we would consider inversion generated by a M\"obius element type,
(\ref{XZ}), say, $F=X_W$ not just $F=X_0$. But then $X_ZX_WX_Z=X_{2Z+W}$, so we can always adsorb $W$ into
$Z$ thus producing no new factors; we therefore stay with our choice of~$F$.},
\begin{gather*}
X_{Z}FX_{Z}=X_{2Z}.
\end{gather*}

We therefore preserve the notion of the geodesic function for curves with inversions as well. We consider
{\em all possible} paths in the spine (graph) that are closed and may experience an arbitrary number
of inversions at pending vertices of the graph. As above, we associate with such paths the geodesic functions
(here, we let $Z_i$ denote the variables of pending edges and $Y_j$ all other variables)
\begin{gather}
\label{Gref}
G_{\gamma}\equiv \tr P_\gamma=2\cosh(\ell_\gamma/2)=\tr LX_{Z_n}FX_{Z_n}R X_{Y_{n-1}}\cdots RX_{Z_1}F X_{Z_1}.
\end{gather}

We have that, for the windowed surface $\Sigma_{g,\delta}$, the number of the shear coordinates $Z_\alpha$ is
\begin{gather*}
\# Z_\alpha=6g-6+3s+2\,\sum_{j=1}^s\delta_j,
\end{gather*}
and adding a new window increases this number by two.

Before describing the general structure of algebras of geodesic functions, let us clarify the geometric
origin of our construction in the simplest possible example.

\subsection{Annulus with one marked point}\label{ss:annul}

\begin{figure}[tb]
\centerline{\includegraphics{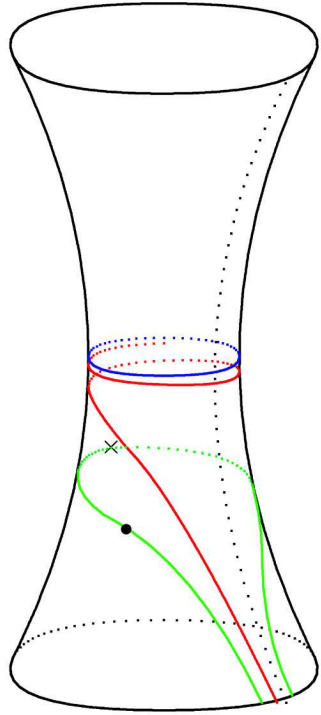}}
\caption{Geodesic lines on the hyperboloid:
dotted vertical line is the asymptote going to the marked point on
the absolute, closed blue line is a unique closed geodesic; red line
is the line from the ideal triangular decomposition asymptotically
approaching the asymptote by one end and the closed geodesic by the
other; green line is the line of inversion whose both ends approach
the marked point; we let the bullet on this line denote the unique
stable point under the inversion and the cross denote the point that
is closest to the closed geodesic.}
\label{fi:hyp-geom}
\end{figure}

The simplest example is the
annulus with one marked point on one of the boundary components (another example of disc
with three marked points will be considered later). Here, the geometry is as in
Fig.~\ref{fi:hyp-geom} where we let the closed line around the neck (the blue line)
denote a unique closed geodesic corresponding
to the element $P_{\mathrm I}$ of the Fuchsian group to be def\/ined below, the winding to it line
(the red line) is the boundary geodesics from the (ideal)
triangle description due to Penner and Fock, and the lower geodesic (the green line) is the new line of inversion (the window).
We indicate by bullet the stable point and by cross the point of the inversion line that is closest to
the closed geodesic.

\begin{figure}[tb]
\centerline{\includegraphics{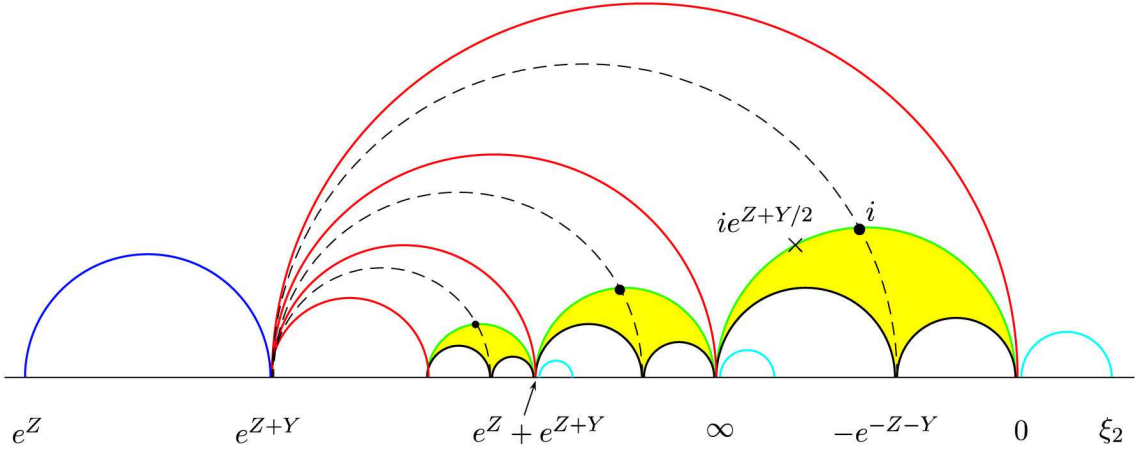}}
\caption{The hyperbolic picture corresponding to the pattern in Fig.~\ref{fi:hyp-geom}:
preimages of red boundary line are red half-circles, preimages of the inversion line
are green half circles (the selected one connects the points $\infty$ and $0$ on the absolute),
and the preimage of the closed geodesic is the (unique) blue half-circle;
the points $e^Z$ and $e^{Z+Y}$ on the absolute are stable under the action of the corresponding Fuchsian
element $P_{\mathrm{I}}$ (\ref{PI}); the bullet symbols are
preimages of the point that is stable upon inversion (the one that lies on the geodesic line
between $\infty$ and $0$ is $i$ in the standard coordinates on the upper half-plane)
and the dotted half-circles connect the point $e^{Z+Y}$
with its images (one of which is
$-e^{-Z-Y}$) under the action of the inversion element $F$. We also mark by cross the
point $ie^{Z+Y/2}$ of the green geodesic line that is closest to the closed geodesic.
The invariant axis of the new element $P_{\mathrm{II}}$ (\ref{PII}) and some of its images under
the action of (\ref{PI}) are depicted as cyan half-circles; $\xi_2$ is from (\ref{xi2}).}
\label{fi:Poincare}
\end{figure}

The same picture in the Poincar\'e upper half-plane is presented  in Fig.~\ref{fi:Poincare}.
There, the whole domain in Fig.~\ref{fi:hyp-geom} bounded below by the bordered (green) geodesic line
and above by the neck geodesic (blue) line is obtained from a single ideal triangle with the vertices
$\{e^{Z+Y}, \infty, 0\}$ upon gluing together two (red) sides of this ideal triangle.
We now construct two (hyperbolic) ele\-ments: $P_{\mathrm I}$ that is the generating element for the original
hyperbolic geometry and the new element
$P_{\mathrm {II}}$ that corresponds to the inversion w.r.t. the lower (green) geodesic
line in Figs.~\ref{fi:hyp-geom} and~\ref{fi:Poincare}.
Adding this new element obviously changes the pattern, but because the Fuchsian pro\-perty retains,
the quotient of the Poincar\'e  upper half-plane under the action of this new Fuchsian group
must be again a Riemann surface with holes. As we demonstrate below, this new Riemann surface is just
the double of the initial bordered Riemann surface.

\begin{figure}[tb]
\centerline{\includegraphics{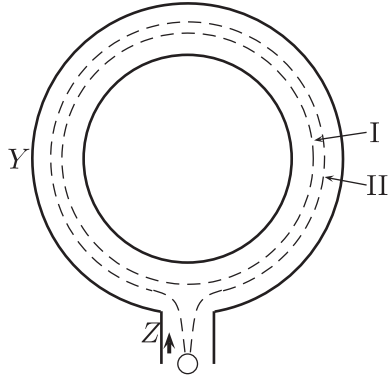}}
\caption{The graph for annulus with one marked point
on one of the boundary components. Examples of closed
geodesics without inversion (I) and with inversion (II) are presented.
The short fat arrow indicates the starting direction for elements of the
Fuchsian group.}
\label{fi:annulus}
\end{figure}

For this, we use the graph representation.
The corresponding fat graph is depicted in Fig.~\ref{fi:annulus}.
This graph with one pending edge and another edge that starts and terminates at the same
vertex is dual to an ideal triangle $\{e^{Z+Y}, \infty, 0\}$ in which two (red) sides are
glued one to another (the resulting curve is dual to the loop) and the remaining (green)
side is the boundary curve (dual to the pending edge).
We mark the starting direction by
the fat arrow, so the element~$P_{\mathrm I}$~is
\begin{gather}
\label{PI}
P_{\mathrm I}=X_ZLX_YLX_Z=
\left(\begin{array}{cc} e^{-Y/2}+e^{Y/2} & -e^{Z+Y/2}\\ e^{-Z-Y/2} & 0\end{array}\right).
\end{gather}
Apparently, the corresponding geodesic function $G_{\mathrm I}$ is just $e^{-Y/2}+e^{Y/2}$, so the
length of the closed geodesic is $|Y|$ as expected.

We now construct the element $P_{\mathrm {II}}$. Note that this element makes
the inversion w.r.t.\ the geodesic between $0$ and $\infty$, so we set the matrix $F$ {\em first}
(since the multiplication is from right to left, this matrix will be rightmost). Then, the rest is
just the above element $P_{\mathrm I}$:
\begin{gather}
\label{PII}
P_{\mathrm {II}}=X_ZLX_YLX_ZF=P_{\mathrm I}F=
\left(\begin{array}{cc} e^{Z+Y/2} & 0\\ e^{-Y/2}+e^{Y/2} & e^{-Z-Y/2}\end{array}\right),
\end{gather}
and the corresponding geodesic function $G_{\mathrm {II}}$ is $2\cosh(Z+Y/2)$ so the length of the
corresponding geodesic (but in a geometry still to be def\/ined!) is $|2Z+Y|$.

\looseness=-1 We now consider the action of these two elements in the geometry of the Poincar\'e upper half-plane
in Fig.~\ref{fi:Poincare}.
It is easy to see that the element $P_{\mathrm I}$ has two stable points: $e^Z$ (attractive) and $e^{Z+Y}$
(repulsive). It also maps $\infty\to 0$, \ $e^Z+e^{Z+Y}\to\infty$, etc. thus producing the
inf\/inite set of preimages of the red geodesic line in Fig.~\ref{fi:hyp-geom} upon identif\/ication under
the action of this element.

The element $F$ f\/irst interchanges $0$ and $\infty$ and $e^{Z+Y}$ and $-e^{-Z-Y}$ thus establishing the
inversion (inversion) w.r.t.\ the green geodesic line. The only stable point of this inversion is the
point of intersection of the two above geodesic lines, and it is the point $i$ in the upper complex half-plane for every
$Z+Y$. Further action is given by $P_{\mathrm I}$ and, in particular, it maps $\infty$ back to~$0$, so $\xi_1=0$
is a stable point of $P_{\mathrm {II}}$. Another stable point is
\begin{gather}
\label{xi2}
\xi_2=\frac{e^{Z+Y/2}-e^{-Z-Y/2}}{e^{Y/2}+e^{-Y/2}},
\end{gather}
and it is easy to see that the two invariant axes of $P_{\mathrm I}$ and $P_{\mathrm {II}}$ never intersect.
Adding the element $P_{\mathrm {II}}$ to the set of generators of the new extended
Fuchsian group we therefore obtain a~new geometry.

\begin{figure}[tb]
\centerline{\includegraphics{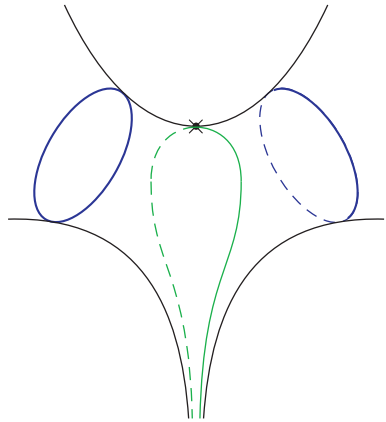}}
\caption{The doubled Riemann surface obtained upon inversion w.r.t.\ the green geodesic
in the case where the stable point coincides with
the point closest to the closed geodesics (blue line)
(the cross then coincides with the bullet).}
\label{fi:symmetry}
\end{figure}

First, let us consider the special case where the stable point on the inversion curve coincides
with the point that is closest to the closed geodesic. Then, apparently, the inversion process
exhibits a symmetry depicted in Fig.~\ref{fi:symmetry}. Considering the Riemann surface
depicted in Fig.~\ref{fi:hyp-geom}, we chop out all its part that is below the green (inversion) line.
We then obtain the {\em double} of the Riemann surface merely by inverting it w.r.t.\ the green line
taking into account the obvious (mirror) symmetry that takes place in this case. We then obtain from
the hyperboloid with marked point at the boundary component the sphere with two identical cycles
(images of the closed geodesic) and one additional puncture (hole of zero length), as shown in
Fig.~\ref{fi:symmetry}.

What happens if, instead of the stable point marked by cross, we have arbitrary stable point (bullet in
Figs.~\ref{fi:hyp-geom} and~\ref{fi:Poincare})? Actually, we can answer this question just from the
geometrical standpoint. Indeed, since in the pattern in Fig.~\ref{fi:hyp-geom}, points on the inversion
geodesics that lie to both sides from the asymptote are close, they must remain close in the new geometry.
But the image of each such point is shifted by a distance that is twice the distance $D$ (along the inversion
line, which is a geodesic line) between the stable point (bullet) and the symmetric point (cross). This means
that, in the new geometry, the points on the inversion line separated by a distance $2D$ must be asymptotically close
as approaching the absolute in the pattern of Fig.~\ref{fi:Poincare}. This means in turn that the corresponding
geodesic in the new geometry is just a geodesic approaching the new closed geodesic of length $2D$.

\begin{figure}[tb]
\centerline{\includegraphics{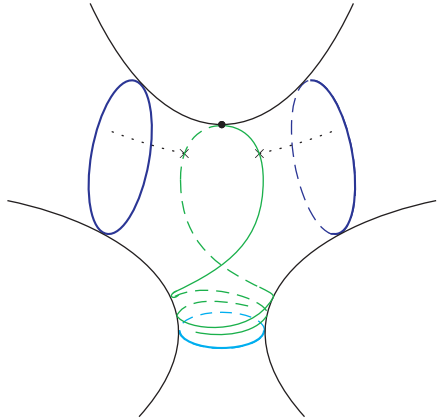}}
\caption{The doubled Riemann surface obtained upon inversion w.r.t.\ the green geodesic
in the case where the stable point (marked by $\bullet$) is arbitrary. The closed
in the asymptotic geodesic sense points in the new geometry are those on dif\/ferent coils
of the spiraling green geodesics, which has the asymptotic form of the double helix.
The separation length is asymptotically equal to $|2Z+Y|$.
The cyan line is the new closed geodesic (the invariant axis of the element $P_{\mathrm{II}}$).
We let two crosses denote the points on the inversion line that are closest to the two copies of
the initial closed geodesic line; the geodesic distance between them is also $|2Z+Y|$.}
\label{fi:tri-sphere}
\end{figure}

It remains just to note that, from the pattern in Fig.~\ref{fi:Poincare},
\[
D=|Z+Y/2|,
\]
that is, the perimeter of the new hole is $|2Z+Y|$, and it coincides with the length of the new element $P_{\mathrm{II}}$
(\ref{PII}), which is therefore the element of the new, extended, Fuchsian group corresponding to going round the new hole.
In Fig.~\ref{fi:tri-sphere}, we depict this new geometry. It is also interesting to note that we now again,
as in the symmetrical case, have two (homeomorphic) images of the initial bordered surface, but the union of these two images
in Fig.~\ref{fi:tri-sphere} constitutes only the part of the corresponding Riemann surface that is above the new
closed geodesics (the cyan line); two ends of the green geodesics constitute the double helix approaching the new geodesic line but
never reaching it, and we always have one copy of the initial surface on one side of coils of this helix and the other copy~-- on
the other side.

\section{Algebras of geodesic functions}\label{s:algebra}

\subsection{Poisson structure}\label{ss:Poisson}

One of the most attractive properties of the graph description is a very simple Poisson algebra on the space
of parameters $Z_\alpha$. Namely, we have the following theorem. It was formulated for surfaces
without marked points in~\cite{Fock1} and here we extend it to {\em arbitrary graphs} with pending
vertices.

\begin{theorem}\label{th-WP} In the coordinates $(Z_\alpha )$ on any fixed spine
corresponding to a surface with marked points on its boundary components,
the Weil--Petersson bracket $B_{{\mbox{\tiny \rm WP}}}$ is given by
\begin{gather}
\label{WP-PB}
B_{{\mbox{\tiny\rm  WP}}}
= \sum_{v} \sum_{i=1}^{3}\frac{\partial}{\partial Z_{v_i}}\wedge
\frac{\partial}{\partial Z_{v_{i+1}}},
\end{gather}
where the sum is taken over all three-valent {\rm(}i.e., not pending\/{\rm)} vertices~$v$ and $v_i$,
$i=1,2,3\ \hbox{\rm mod}\ 3$, are the labels of the cyclically ordered
edges incident on this vertex irrespectively on whether they are internal or
pending edges of the graph.
\end{theorem}

The center of this Poisson algebra is provided by the proposition.

\begin{proposition}\label{prop12}
The center of the Poisson algebra \eqref{WP-PB}) is generated by
elements of the form $\sum Z_\alpha$, where the sum is over all edges
of $\Gamma $ in a boundary component of $F(\Gamma )$ taken with multiplicities.
This means, in particular, that each pending edge contributes twice to such sums.
\end{proposition}

\begin{proof} The proof is purely technical; for the case of surfaces without marked
points on boundary components it can be found in Appendix~B in \cite{ChP}. When adding
marked points, it is straightforward to verify that the sums in the assertion of the
proposition are central elements. In order to prove that no extra central elements
appear due to the addition process, it suf\/f\/ices to verify that the two changes of
the part of a graph shown below,

\centerline{\includegraphics{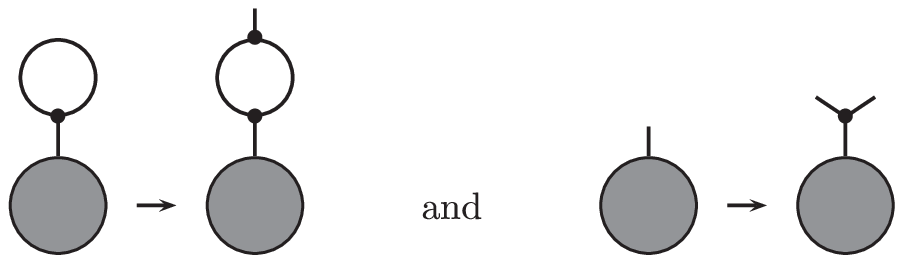}}

\noindent
do not change the corank of the Poisson relation matrix $B(\Gamma_{g,s,\delta})$.\end{proof}

\begin{example}\label{Ex1}
Let us consider the graph in Fig.~\ref{fi:center}. It has two boundary components and two
corresponding geodesic lines. Their lengths, $\sum\limits_{i=1}^4 Y_i$ and $\sum\limits_{i=1}^4 (Y_i+2Z_i)$,
are the two Casimirs of the Poisson algebra with the def\/ining relations
\[
\{Y_i,Y_{i-1}\}=1\quad\mod\ 4,\qquad \{Z_i,Y_i\}=-\{Z_i,Y_{i-1}\}=1\quad \mod 4,
\]
and with all other brackets equal to zero.
\end{example}

\subsection[Classical flip morphisms and invariants]{Classical f\/lip morphisms and invariants}\label{ss:flip}

The $Z_\alpha$-coordinates (which are the logarithms of cross ratios) are called {\it (Thurston) shear
coordinates} \cite{ThSh,Bon2} in the case of punctured Riemann surface (without boundary components).
We preserve this notation and this term also in the case of windowed surfaces.

In the case of surfaces with holes, $Z_\alpha$ were the coordinates on the Teichm\"uller space
${\mathcal T}_{g,s}^H$, which was the $2^s$-fold covering of the standard Teichm\"uller space ramif\/ied over
surfaces with punctures (when a hole perimeter becomes zero, see~\cite{Fock2}). We assume correspondingly
$Z_\alpha$ to be the coordinates of the corresponding spaces ${\mathcal T}_{g,\delta}^H$ in the bordered
surfaces case.

Assume that there is an enumeration of the edges of $\Gamma $ and that
edge $\alpha$ has distinct endpoints. Given a spine $\Gamma$ of $\Sigma$, we may produce
another spine $\Gamma _\alpha$ of $\Sigma$ by contracting and expanding edge~$\alpha$ of
$\Gamma $, the edge labelled $Z$ in Fig.~\ref{fi:flip}, to produce $\Gamma _\alpha$ as in the
f\/igure; the fattening and embedding of $\Gamma_\alpha$ in $\Sigma$ is determined from
that of $\Gamma$ in the natural way. Furthermore, an enumeration of the edges of
$\Gamma $ induces an enumeration of the edges of $\Gamma _\alpha$ in the natural
way, where the vertical edge labelled $Z$ in Fig.~\ref{fi:flip} corresponds to the horizontal
edge labelled $-Z$.  We say that $\Gamma _\alpha$ arises from~$\Gamma$ by a
{\it Whitehead move} along edge $\alpha$.  We also write $\Gamma
_{\alpha\beta}=(\Gamma _\alpha )_\beta$, for any two indices~$\alpha$,~$\beta$ of
edges, to denote the result of f\/irst performing a move along $\alpha$ and then along
$\beta$; in particular, $\Gamma _{\alpha \alpha}=\Gamma$ for any index $\alpha$.

\subsubsection{Whitehead moves on inner edges}\label{sss:mcg}

\begin{figure}[tb]
\centerline{\includegraphics{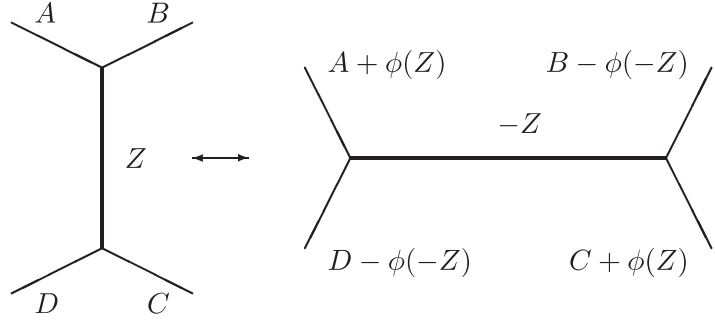}}
\caption{Flip, or Whitehead move on the shear coordinates $Z_\alpha$. The outer edges can be
pending, but the inner edge with respect to which the morphism is performed cannot be a pending
edge.}
\label{fi:flip}
\end{figure}

\begin{proposition}[\cite{ChF}]\label{propcase}
Setting $\phi (Z)={\rm log}(e^Z+1)$ and adopting the notation of Fig.~{\rm \ref{fi:flip}}
for shear coordinates of nearby edges, the effect of a
Whitehead move is as follows:
\begin{gather}
W_Z\,:\ (A,B,C,D,Z)\to (A+\phi(Z), B-\phi(-Z), C+\phi(Z), D-\phi(-Z), -Z).
\label{abc}
\end{gather}
In the various cases where the edges are not distinct
and identifying an edge with its shear coordinate in the obvious notation we have:
if $A=C$, then $A'=A+2\phi(Z)$;
if $B=D$, then $B'=B-2\phi(-Z)$;
if $A=B$ (or $C=D$), then $A'=A+Z$ (or $C'=C+Z$);
if $A=D$ (or $B=C$), then $A'=A+Z$ (or $B'=B+Z$).
Any variety of edges among $A$, $B$, $C$, and $D$ can be pending edges of the graph.
\end{proposition}

We also have two simple but important lemmas establishing the properties of
invariance w.r.t.\ the f\/lip morphisms.

\begin{lemma} \label{lem1}
Transformation~\eqref{abc} preserves
the traces of products over paths \eqref{Gref}.
\end{lemma}

\begin{lemma} \label{lem-Poisson}
Transformation~\eqref{abc} preserves
Poisson structure \eqref{WP-PB} on the shear coordinates.
\end{lemma}

That the Poisson algebra for the bordered surfaces case
is invariant under the f\/lip transformations follows immediately
because we f\/lip here inner, not pending, edges of a graph, which reduces
the situation to the ``old" statement for surfaces without windows.

We also have the statement concerning the polynomiality of geodesic functions.

\begin{proposition}\label{prop-Laurent}
All $G_{\gamma}$ constructed by \eqref{Gref}
are Laurent polynomials in $e^{Z_i}$ and $e^{Y_j/2}$ with positive integer coefficients, that is, we have
the Laurent property, which holds, e.g., in cluster algebras~{\rm \cite{FZ}}.
All these geodesic functions preserve their
polynomial structures upon Whitehead moves on inner edges, and all of them are hyperbolic elements
$(G_\gamma>2)$, the only exception where $G_\gamma=2$ are paths homeomorphic to going around holes of zero length
(punctures).
\end{proposition}

\begin{figure}[tb]
\centerline{\includegraphics{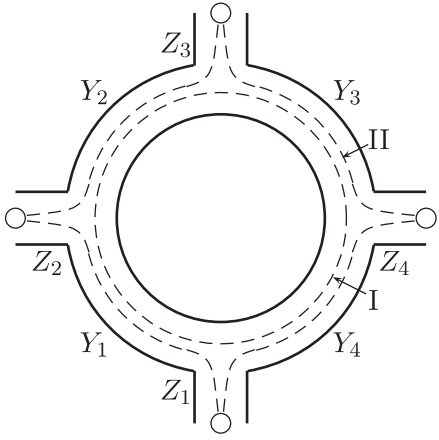}}
\caption{An example of geodesics whose geodesic functions are in
the center of the Poisson algebra (dashed lines). Whereas $G_{\mathrm I}$
corresponds to the standard geodesic around the hole (no marked points are present
on the corresponding boundary component), the line that is parallel
to a boundary component with marked points must experience
all possible inversions on its way around the boundary component, as is the
case for $G_{\mathrm {II}}$.}
\label{fi:center}
\end{figure}

\subsubsection{Whitehead moves on pending edges}\label{sss:pending}

In the case of windowed surfaces, we encounter a new phenomenon as compared with the case
of surfaces with holes. We can construct morphisms relating {\em any} of the Teichm\"uller spaces
${\mathcal T}_{g,\delta^1}^H$ and ${\mathcal T}_{g,\delta^2}^H$ with  $\delta^1=\{\delta_1^1,\dots, \delta_{n_1}^1\}$
and $\delta^2=\{\delta_1^2,\dots, \delta_{n_2}^2\}$ providing $n_1=n_2=n$ and
$\sum\limits_{i=1}^{n_1}\delta^1_i=\sum\limits_{i=1}^{n_2}\delta^2_i$, that is, we explicitly construct morphisms relating
any two of algebras corresponding to windowed surfaces of the same genus, same number of boundary components, and
with the same total number of windows; the window distribution into the boundary components can be however
arbitrary.

This new morphism corresponds in a sense to f\/lipping a pending edge.

\begin{figure}[tb]
\centerline{\includegraphics{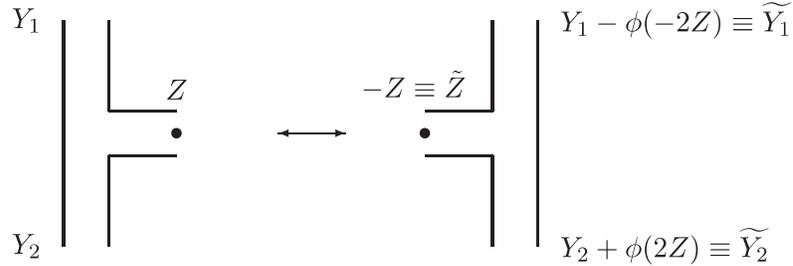}}
\caption{Flip, or Whitehead move on the shear coordinates when f\/lipping the pending
edge $Z$ (indicated by bullet). Any (or both) of edges $Y_1$ and $Y_2$ can be pending.}
\label{fi:mcg-pending}
\end{figure}

\begin{lemma} \label{lem-pending}
Transformation~in Fig.~{\rm \ref{fi:mcg-pending}} is the morphism between the spaces
${\mathcal T}_{g,\delta^1}^H$ and ${\mathcal T}_{g,\delta^2}^H$. These
morphisms preserve both Poisson structures \eqref{WP-PB} and the geodesic
length functions. In Fig.~{\rm \ref{fi:mcg-pending}} any (or both) of $Y$-variables can be
variables of pending edges (the transformation formula is insensitive to it).
\end{lemma}

\begin{proof} Verifying the preservation of Poisson relations (\ref{WP-PB}) is simple, whereas
for traces over paths we have four dif\/ferent cases
of path positions in the subgraph in the left side of Fig.~\ref{fi:mcg-pending}, and in each case
we have the corresponding path in the right side of this f\/igure\footnote{We can think about the
f\/lip in Fig.~\ref{fi:mcg-pending} as about ``rolling the bowl" (the dot-vertex) from one side to the
other; the pending edge is then ``plumbed" on the left and is protruded from the right side whereas
threads of all geodesic lines are deformed continuously, see the example in Fig.~\ref{fi:fol-pend}.}.
In each of these cases we have the following
{\em matrix} equalities (each can be verif\/ied directly)
\begin{gather*}
X_{Y_2}LX_ZFX_ZLX_{Y_1}=X_{{\tilde Y}_2}LX_{{\tilde Y}_1},\nonumber\\
X_{Y_1}RX_ZFX_ZRX_{Y_1}=X_{{\tilde Y}_1}LX_{\tilde Z}FX_{\tilde Z}RX_{{\tilde Y}_1},\nonumber\\
X_{Y_2}RX_{Y_1}=X_{{\tilde Y}_2}RX_{\tilde Z}FX_{\tilde Z}RX_{{\tilde Y}_1},\nonumber\\
X_{Y_2}LX_ZFX_ZRX_{Y_2}=X_{{\tilde Y}_2}RX_{\tilde Z}FX_{\tilde Z}LX_{{\tilde Y}_2},\nonumber
\end{gather*}
where (in the exponentiated form)
\begin{gather*}
e^{{\tilde Y}_1}=e^{Y_1}\bigl(1+e^{-2Z}\bigr)^{-1},\qquad e^{{\tilde Y}_2}=e^{Y_1}\bigl(1+e^{2Z}\bigr),\qquad
e^{{\tilde Z}}=e^{-Z}.\tag*{\qed}
\end{gather*}\renewcommand{\qed}{}
\end{proof}

From the technical standpoint, all these equalities follow from f\/lip transformation (\ref{abc})
upon the substitution $A=C=Y_2$, $B=D=Y_1$, and $Z=2Z$. The above four cases of geodesic functions are
then exactly four possible cases of geodesic arrangement in the (omitted) proof of Lemma~\ref{lem1}.

Using f\/lip morphisms in Fig.~\ref{fi:mcg-pending} and in formula (\ref{abc}), we may establish a morphism
between any two algebras corresponding to surfaces of the same genus, same number of boundary components, and
same total number of marked points on these components; their distribution into the boundary components can be
however arbitrary. And it is again a standard tool that if, after a series of morphisms, we come to a graph
of the same combinatorial type as the initial one (disregarding marking of edges), we associate a {\em mapping class group}
operation with this morphism therefore passing from the groupoid of morphisms to the group of
modular transformations.

\begin{example}\label{Ex2}
The f\/lip morphism w.r.t.\ the edge $Z_1$ in the pattern in (\ref{braid1}),
\begin{gather}\label{braid1}
\begin{gathered}[c]
\includegraphics{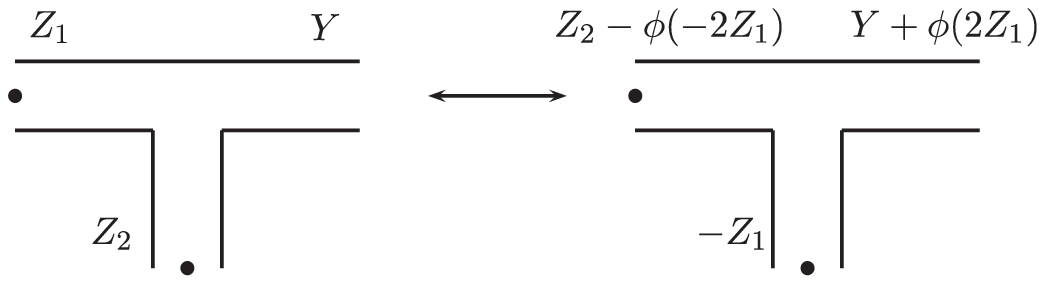}
\end{gathered}
\end{gather}
where $Z_1$ and $Z_2$ are the pending edges, generates the (unitary) mapping class group transformation
\begin{gather*}
e^{Z_2}\to e^{-Z_1},\qquad e^{Z_1}\to e^{Z_2}\bigl(1+e^{-2Z_1}\bigr)^{-1},\qquad e^Y\to e^Y\bigl(1+e^{2Z_1}\bigr)
\end{gather*}
on the corresponding Teichm\"uller space ${\mathcal T}_{g,\delta}^H$. This is a particular case of {\em braid}
transformation to be considered in detail in Section~\ref{ss:braid}.
\end{example}

\subsection{New graphical representation}\label{ss:new}

\begin{figure}[tb]
\centerline{\includegraphics{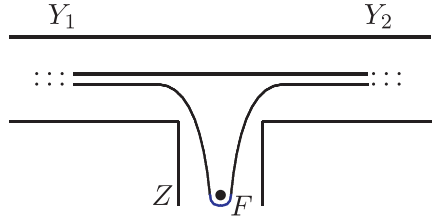}}
\caption{Resolution of the inversion process from Fig.~\ref{fi:corner}.
We introduce the new dot-vertex reducing the inversion to winding around
this vertex (blue part of the path in the graph). Each time a path winds
around the dot-vertex, we set the inversion matrix $F$.}
\label{fi:dot-vertex}
\end{figure}

In the case of usual geodesic functions, there exists a very convenient representation in which
one can apply classical skein and Poisson relations in classical case or the quantum skein
relation in the quantum case and ensure the Riedemeister moves when ``disentangling" the products of
geodesic function representing them as linear combinations of multicurve functions.
However, in our case, it is still obscure what happens when geodesic lines
intersect in some way at the pending vertex. In fact, we can propose the comprehensive graphical
representation in this case as well! For this, let us come back to Fig.~\ref{fi:corner} and resolve now
the inversion introducing a new {\em dot-vertex} at a pending vertex inside the fat graph
and assuming that the inversion
matrix $F$ corresponds to winding around this dot-vertex as shown in Fig.~\ref{fi:dot-vertex}.

We now formulate the rules for geodesic algebra that follow from relations (\ref{WP-PB}) and
classical skein relations. They coincide with the rules in the case of surfaces with holes except
the one new case depicted in Fig.~\ref{fi:dot-skein}. Note that all claims below follow from
direct and explicit calculations involving representations from Section~\ref{s:hyp}.

\subsubsection{Classical skein relation}\label{sss:skein}

The trace relation $\tr(AB)+\tr(AB^{-1})-\tr A\cdot\tr B=0$ for
arbitrary $2\times 2$ matrices~$A$ and~$B$ with unit determinant allows one to
``disentangle'' any product of geodesic functions,
i.e., express it uniquely as a f\/inite linear combination of generalized multicurves
(see Def\/inition~\ref{def1} below).
Introducing the additional factor $\#G$ to be the total number of components
in a multicurve, we can uniformly present the classical skein relation as
\begin{gather}
\begin{gathered}[c]
\includegraphics{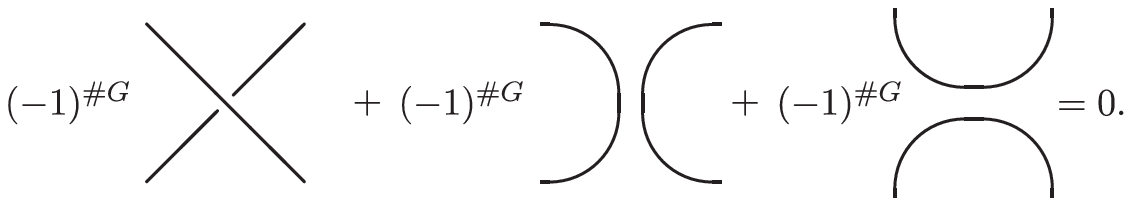}
\end{gathered}
\label{skeinclass}
\end{gather}

We assume in (\ref{skeinclass}) that the ends of lines are joint pairwise in the rest of the
graph, which is the same for all three items in the formula. Of course, we perform
there algebraic operations with the algebraic quantities~-- with the (products of) geodesic
functions corresponding to the respective families of curves.

\subsubsection{Poisson brackets for geodesic functions}

We f\/irst mention that two geodesic functions Poisson
commute if the underlying geodesics are disjointly embedded in the sense of the new
graph technique involving dot-vertices. Because of
the Leibnitz rule for the Poisson bracket, it suf\/f\/ices to consider only ``simple'' intersections of
pairs of geodesics with respective geodesic functions
$G_1$ and $G_2$ of the form
\begin{gather}
\label{G1}
G_1=\tr^1\cdots X_C^1R^1X_Z^1L^1X_A^1\cdots,
\\
\label{G2}
G_2=\tr^2\cdots X_B^2 L^2 X_Z^2R^2X_D^2\cdots,
\end{gather}
where the superscripts $1$ and~$2$ pertain to operators and traces in two dif\/ferent
matrix spaces.

The positions of edges $A,B,C,D,$ and $Z$ are as in Fig.~\ref{fi:flip}.
Dots in (\ref{G1}), (\ref{G2}) refer to arbitrary sequences of matrices $R^{1,2}$, $L^{1,2}$,
$X_{Z_i}^{1,2}$, and $F^{1,2}$ belonging to the corresponding matrix spaces; $G_1$ and $G_2$ must correspond
to closed geodesic lines, but we make no assumption on their simplicity or graph simplicity, that is,
the paths that correspond to $G_1$ and $G_2$ may have self- and mutual intersections and,
in particular, may pass arbitrarily many times through the edge $Z$ in Fig.~\ref{fi:flip}.

Direct calculations then give
\begin{gather}
\label{Goldman}
\{G_1,G_2\}=\frac{1}{2}(G_{\mbox{\tiny H}}-G_{\mbox{\tiny I}}),
\end{gather}
where $G_{\mbox{\tiny I}}$ corresponds to the geodesic that is obtained by
erasing the edge~$Z$ and joining together the edges ``$A$'' and ``$D$'' as
well as ``$B$'' and ``$C$'' in a natural way as illustrated in the middle diagram
in (\ref{skeinclass}); $G_{\mbox{\tiny H}}$ corresponds to the geodesic that
passes over the edge~$Z$ twice, so it has the form
$\tr \cdots X_CR_ZR_D\cdots$ $\cdots X_BL_ZL_A\cdots$
as illustrated in the rightmost diagram in (\ref{skeinclass}).
These relations were f\/irst obtained
in~\cite{Gold} in the continuous parametrization (the classical Turaev--Viro
algebra).

\begin{figure}[tb]
\centerline{\includegraphics{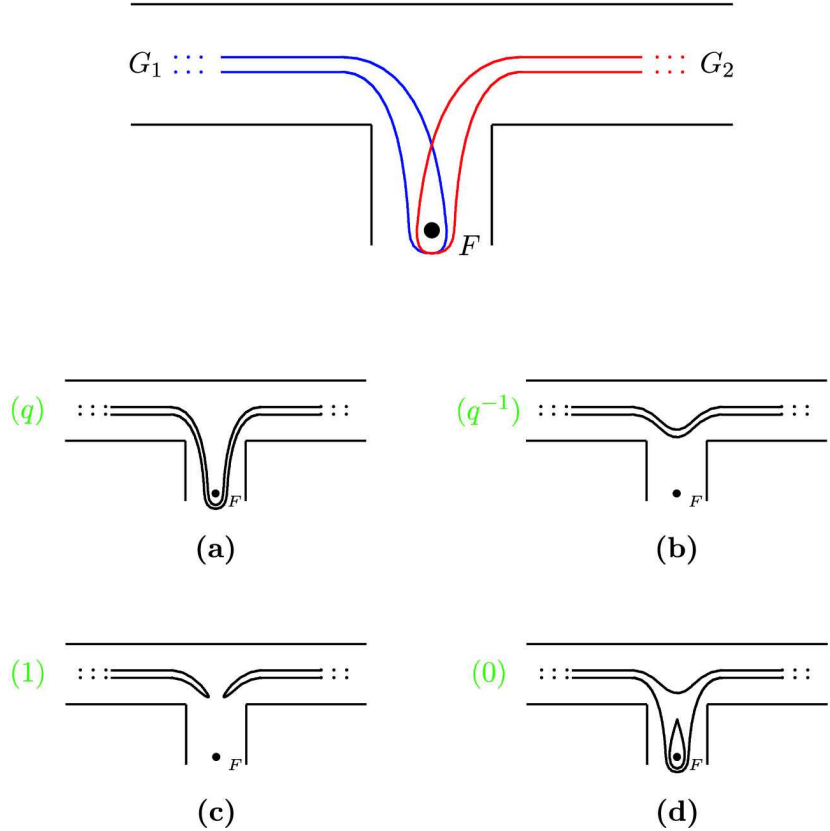}}
\caption{An example of two geodesic lines intersecting at the
dot-vertex. We present four homotopical types of resolving two
intersections in this pattern (Cases (a)--(d)). Case (d) contains the loop with only the dot-vertex inside.
This loop is $\tr F=0$, so the whole contribution vanishes in this case. The (green) factors
in brackets pertain to the quantum case in Section~\ref{ss:q} indicating the weights with
which the corresponding (quantum) geodesic multicurves enter the expression for the
product $G_1^\hbar G_2^\hbar$.}
\label{fi:dot-skein}
\end{figure}

Having two curves, $\gamma_1$ and $\gamma_2$, with an arbitrary number of crossings,
we now f\/ind their Poisson bracket using the following rules:
\begin{itemize}\itemsep=0pt
\item We take a sum of products of geodesic functions
of non(self)intersecting curves obtained when we apply Poisson relation~(\ref{Goldman})
at one intersection point and classical skein relation~(\ref{skeinclass})
at all the remaining points of intersection; we assume the summation over
all possible cases.
\item If, in the course of calculation, we meet an empty (contractible) loop, then
we associate the factor $-2$ to such a loop; this assignment, as is known~\cite{ChP},
ensures the Riedemeister moves on the set of geodesic lines thus making the bracket
to depend only on the homotopical class of the curve embedding in the surface.
\item If, in the course of calculation, we meet a curve homeomorphic to passing around a
dot-vertex, then we set $\tr F=0$ into the correspondence to such curve thus killing the whole corresponding
multicurve function.
\end{itemize}

These simple and explicit rules are an ef\/fective tool for
calculating the Poisson brackets in many important cases below.

Because the Poisson relations are completely determined by homotopy types of curves involved,
using Lemma~\ref{lem-pending}, we immediately come to the following theorem

\begin{theorem}\label{th-mark}
Poisson algebras of geodesic functions for the bordered Riemann surfaces $\Sigma_{g,\delta^1}$
and $\Sigma_{g,\delta^2}$ that differ only by distributions of marked points among their boundary
components are isomorphic; the isomorphism is described by Lemma~{\rm \ref{lem-pending}}.
\end{theorem}

It follows from this theorem that we can always collect all the marked points on just one boundary
component.

\subsection[The $A_n$ algebras]{The $\boldsymbol{A_n}$ algebras}\label{ss:An}

\begin{figure}[tb]
\centerline{\includegraphics{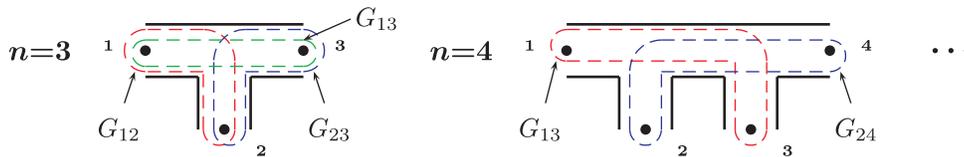}}
\caption{Generating graphs for $A_n$ algebras for $n=3,4,\dots$. We indicate character
geodesics whose geodesic functions $G_{ij}$ enter bases of the corresponding algebras.}
\label{fi:An}
\end{figure}

Consider the disc with $n$ marked points on the boundary; examples of the corresponding
representing graph $\Gamma_n$ are depicted in Fig.~\ref{fi:An} for $n=3,4,\dots$.
We enumerate the $n$ dot-vertices clockwise, $i,j=1,\dots, n$.
We then let $G_{ij}$ with $i<j$ denote the geodesic function corresponding to the geodesic
line that encircles exactly two dot-vertices
with the indices $i$ and $j$. Examples are in the f\/igure: for $n=3$, red line corresponds
to $G_{12}$, blue~-- to $G_{23}$ and green~-- to $G_{13}$. Note that in the cluster terminology
(see \cite{FST}) these algebras were called the $A_{n-2}$-algebras.

Using the skein relation, we can
close the Poisson algebra thus obtaining for $A_3$:
\begin{gather}
\label{alg-A3}
\left\{G_{12},G_{23}\right\}=G_{12}G_{23}-2G_{13}\hbox{ and cycl. permut.}
\end{gather}
Note that the left-hand side is doubled in this case as compared to Nelson--Regge algebras
recalled in~\cite{ChP}. In the $A_3$, case this is easily understandable because, say,
\begin{gather}
\label{G12}
G_{12}=\tr LX_{2Z_2}RX_{2Z_1}=e^{Z_1+Z_2}+e^{Z_1-Z_2}+e^{-Z_1-Z_2},
\end{gather}
and this expression literally coincides with the one
for the algebra of geodesics in the case of higher genus surfaces with one or two holes
(see~\cite{ChF2}) but the left-hand side of the relation is now doubled (the analogous
expression for $G_{12}$ in \cite{ChF2} was the same as in (\ref{G12}) upon the substitution
$Z_1=X_1/2$ and $Z_2=X_2/2$, but with the $X$-variables having the doubled Poisson brackets
$\{X_2,X_1\}=2$).
In higher-order algebras (starting with $n=4$), we meet
a more complicate case of the fourth-order crossing (as shown in the case $n=4$ in Fig.~\ref{fi:An}).
Using our rules for Poisson brackets, we f\/ind that those for these geodesic functions are
\begin{gather}
\label{alg-A4}
\left\{G_{13},G_{24}\right\}=2G_{12}G_{34}-2G_{14}G_{23}
\end{gather}
(note that the items in the products in the r.h.s.\ mutually commute).

It is also worth mentioning that after this doubling that occurs in the right-hand sides of relations
(\ref{alg-A3}) and (\ref{alg-A4}), we come exactly to algebras appearing in the Frobenius manifold
approach~\cite{DM}.

\subsection[The $D_n$-algebras]{The $\boldsymbol{D_n}$-algebras}

\begin{figure}[tb]
\centerline{\includegraphics{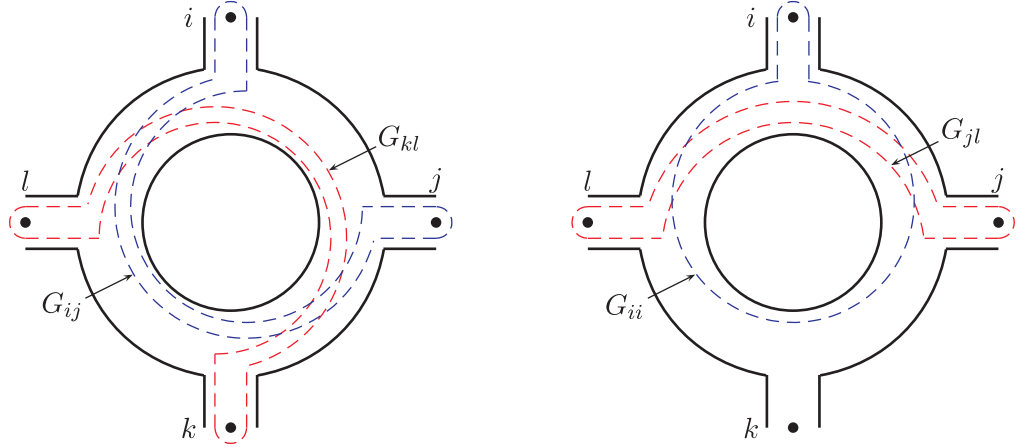}}
\caption{Typical geodesics corresponding to the geodesic functions
constituting a set of generators of the $D_n$ algebra. We let $G_{ij}$,
$i,j=1,\dots,n$, denote these functions. The order of subscripts is now important:
it indicates the direction of encompassing the hole (the second boundary component of
the annulus). The most involved pattern of intersection is on the right part of the f\/igure:
the geodesics have there eight-fold intersection; in the left part we present also the
geodesic function $G_{ii}$ corresponding to the geodesic that starts and terminates at the same
window.}
\label{fi:Dn}
\end{figure}

We now consider the case of annulus with $n$ marked points on one of the boundary component
(see the example in Fig.~\ref{fi:center}. Here, again, the state of art is to f\/ind a convenient (f\/inite)
set of geodesic functions closed w.r.t.\ the Poisson brackets\footnote{Usually we can say nothing
about the uniqueness of such a set for a particular geometry.}. In the case of annulus, such a set
is given by geodesic functions corresponding to geodesics in Fig.~\ref{fi:Dn}.

We therefore describe a set of geodesic functions by the matrix $G_{ij}$ with $i,j=1,\dots, n$ where
the order of indices indicates the direction of encompassing the second boundary component of the annulus.

\begin{lemma} \label{lem-Dn-Poisson}\looseness=-1
The set of geodesic functions $G_{ij}$ corresponding to geodesics in Fig.~{\rm \ref{fi:Dn}} is
Poisson closed.
\end{lemma}

The relevant Poisson brackets are too cumbersome and we omit them here because one can
easily read them from the corresponding quantum algebra in formula~(\ref{Dn-q})
below in the limit as $\hbar\to0$.

\subsection{Braid group relations for windowed surfaces}\label{ss:braid}

\subsubsection[Braid group relations on the level of $Z$-variables]{Braid group relations on the level of $\boldsymbol{Z}$-variables}

\begin{figure}[tb]
\centerline{\includegraphics{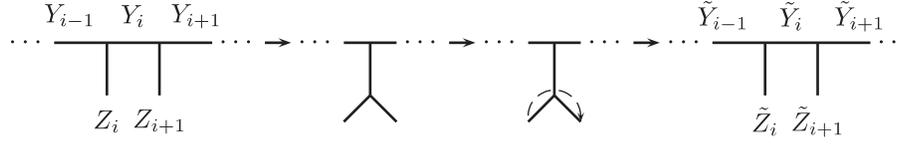}}
\caption{Three-step f\/lip transformation of intertwining pending edge
variables $Z_i$ and $Z_{i+1}$ that results in the same combinatorial
graph. The rest of the graph denoted by dots remains unchanged.}
\label{fi:3step}
\end{figure}

We have already demonstrated in Example \ref{Ex2} a m.c.g.\ relation interchanging two pending
edges of a graph. In a more general case of $A_n$-algebra, we have a graph depicted in Fig.~\ref{fi:An}
and another intertwining relation arises from the three-step f\/lipping process schematically depicted in
Fig.~\ref{fi:3step}.

The graph for the $A_n$ algebra has the form in Fig.~\ref{fi:An} with $Y_i$,
$2\le i\le n-2$, being the variables of internal edges and $Z_j$, $1\le j\le n$,
being the variables of the pending edges and we identify $Y_1\equiv Z_1$ and
$Y_{n-1}\equiv Z_n$ to make formulas below uniform.

We let $R_{i,i+1}$ denote the intertwining transformation in Fig.~\ref{fi:3step} for
$2\le i\le n-2$ and in Fig.~\ref{fi:mcg-pending} for $i=1$ and $i=n-1$. For the
exponentiated variables, these transformations have the form
\begin{gather}
\label{Rii+1}
R_{i,i+1}\left\{
         \begin{array}{l}
           e^{Y_{i-1}} \\
           e^{Y_i}\\
           e^{Y_{i+1}}\\
           e^{Z_i}\\
           e^{Z_{i+1}}
         \end{array}
         \right\}=
\left\{
         \begin{array}{l}
           e^{Y_{i-1}}\bigl(1+e^{2Z_i}(1+e^{Y_i})\bigr) \\
           e^{Y_i}\bigl(1+e^{2Z_i}(1+e^{Y_i})^2\bigr)^{-1}\\
           e^{Y_{i+1}}{1+e^{2Z_i}(1+e^{Y_i})^2\over 1+e^{2Z_i}(1+e^{Y_i})}\\
           e^{2Z_i+Z_{i+1}+Y_i}\bigl(1+e^{2Z_i}(1+e^{Y_i})\bigr)^{-1}\\
           e^{-Z_{i}-Y_i}\bigl(1+e^{2Z_i}(1+e^{Y_i})\bigr)
         \end{array}
         \right\},\qquad 2\le i\le n-2
\end{gather}
and
\begin{gather}
\label{R-12}
R_{1,2}\left\{
         \begin{array}{l}
           e^{Z_{1}} \\
           e^{Z_2}\\
           e^{Y_{2}}
         \end{array}
         \right\}=
\left\{
         \begin{array}{l}
           e^{Z_{2}}(1+e^{-2Z_1})^{-1} \\
           e^{-Z_1}\\
           e^{Y_{2}}(1+e^{2Z_1})
         \end{array}
         \right\},\\
\label{R-n-1n}
R_{n-1,n}\left\{
         \begin{array}{l}
           e^{Z_{n-1}} \\
           e^{Z_n}\\
           e^{Y_{n-2}}
         \end{array}
         \right\}=
\left\{
         \begin{array}{l}
           e^{Z_{n}}(1+e^{-2Z_{n-1}})^{-1} \\
           e^{-Z_{n-1}}\\
           e^{Y_{n-2}}(1+e^{2Z_{n-1}})
         \end{array}
         \right\}.
\end{gather}
The following lemma is the direct calculation using (\ref{Rii+1}), (\ref{R-12}), and (\ref{R-n-1n}).

\begin{lemma} \label{lem-braid}
For any $n\ge3$, we have the {\em braid group relation}
\begin{gather*}
R_{i-1,i}R_{i,i+1}R_{i,i-1}=R_{i,i+1}R_{i-1,i}R_{i,i+1},\qquad 2\le i\le n-1.
\end{gather*}
\end{lemma}

\subsubsection{Braid group relations for geodesic functions of $A_n$-algebras}

Here we, following Bondal~\cite{Bondal}, propose another, simpler way to derive
the braid group relations using the construction of the groupoid of upper-triangular matrices.
It was probably f\/irst used in \cite{DM} to prove the braid group relations in the case of $A_3$ algebra.
In the case of $A_n$ algebras for general $n$, let us construct the upper-triangular matrix ${\mathcal A}$
\begin{gather}
\label{A-matrix}
{\mathcal A}=\left(\begin{array}{ccccc}
                     1 & G_{1,2} & G_{1,3} & \dots & G_{1,n} \\
                     0 & 1 & G_{2,3} & \dots & G_{2,n} \\
                     0 & 0 & 1 & \ddots & \vdots \\
                     \vdots & \vdots & \ddots & \ddots & G_{n-1,n} \\
                     0 & 0 & \dots & 0 & 1 \\
                   \end{array}
\right)
\end{gather}
associating the entries $G_{i,j}$ with the geodesic functions. Using the skein relation,
we can then present the action of the braid group element $R_{i,i+1}$ exclusively in terms
of the geodesic functions from this, f\/ixed, set:
\begin{gather*}
R_{i,i+1}{\mathcal A}={\tilde{\mathcal A}},\qquad \hbox{where} \ \
\left\{
\begin{array}{ll}
  {\tilde G}_{i+1,j}=G_{i,j}, & j>i+1,\\
  {\tilde G}_{j,i+1}=G_{j,i}, & j<i, \\
  {\tilde G}_{i,j}=G_{i,j}G_{i,i+1}-G_{i+1,j}, & j>i+1, \\
  {\tilde G}_{j,i}=G_{j,i}G_{i,i+1}-G_{j,i+1}, & j<i, \\
  {\tilde G}_{i,i+1}=G_{i,i+1}. &  \\
\end{array}%
\right.
\end{gather*}
A very convenient way to present this transformation is by introducing the special matrices
$B_{i,i+1}$ of the block-diagonal form
\begin{gather*}
B_{i,i+1}=\begin{array}{c}
            \vdots \\
            i \\
            i+1 \\
            \vdots \\
          \end{array}
          \left(
          \begin{array}{cccccccc}
            1 &  &  &  &  &  &  &  \\
             & \ddots &  &  &  &  &  &  \\
             &  & 1 &  &  &  &  &  \\
             &  &  &   G_{i,i+1} & -1 &  & &  \\
            &  &  & 1 & 0 &  &  &  \\
             &  &  &  &  & 1 &  &  \\
             &  &  &  &  &  & \ddots &  \\
             &  &  &  &  &  &  & 1 \\
          \end{array}
          \right).
\end{gather*}
Then, The action of the braid group generator $R_{i,i+1}$ on ${\mathcal A}$ is merely
\begin{gather}
\label{BAB}
R_{i,i+1}{\mathcal A}=B_{i,i+1}{\mathcal A}B^{T}_{i,i+1}
\end{gather}
with $B^{T}_{i,i+1}$ the matrix transposed to $B_{i,i+1}$. The proof~\cite{Bondal} of Lemma~\ref{lem-braid} in this setting
is much simpler than in terms of the Teichm\"uller space variables; moreover, using this approach, we can attack another
important issue related to the second braid group relation\footnote{This relation was not presented explicitly in~\cite{Bondal},
so we consider it here in more details.}.

We consider the action of the chain of transformations $R_{n-1,n}R_{n-2,n-1}\cdots R_{2,3}R_{1,2}{\mathcal A}$.
Note that, on each step, the item $G^{(i-1)}_{i,i+1}$ in the corresponding matrix $B_{i,i+1}$ is the transformed quantity
(we assume $G^{(0)}_{ij}$ to coincide with the initial $G_{ij}$ in ${\mathcal A}$).
However, it is easy to see that for just this chain of transformations, $G^{(i-1)}_{i,i+1}=G^{(0)}_{1,i+1}=G_{1,i+1}$,
and the whole chain of matrices $B$ can be then expressed in terms of the initial variables $G_{i,j}$ as
\begin{gather*}
{\mathcal B}\equiv B_{n-1,n}B_{n-2,n-1}\cdots B_{2,3}B_{1,2}=\left(%
\begin{array}{ccccc}
  G_{1,2} & -1 & 0 & \dots & 0 \\
  G_{1,3} & 0 & -1 &  & \vdots \\
  \vdots & \vdots & \ddots & \ddots & 0 \\
  G_{1,n} & 0 & \dots & 0 & -1 \\
  1 & 0 & \dots & 0 & 0 \\
\end{array}%
\right),
\end{gather*}
and the whole action on ${\mathcal A}$ gives
\begin{gather*}
{\tilde {\mathcal A}}\equiv {\mathcal B}{\mathcal A}{\mathcal B}^{T}=\left(%
\begin{array}{cccccc}
  1 & G_{2,3} & G_{2,4} & \dots & G_{2,n} & G_{1,2} \\
  0 & 1 & G_{3,4} & \dots & G_{3,n} & G_{1,3} \\
  0 & 0 & 1 &  & G_{4,n} & G_{1,4} \\
  \vdots & \vdots &  & \ddots &  & \vdots \\
  0 & 0 & \dots & 0 & 1 & G_{1,n} \\
  0 & 0 & \dots & \dots & 0 & 1 \\
\end{array}%
\right),
\end{gather*}
and we see that it boils down to a mere permutation of the elements of the initial matrix ${\mathcal A}$.
It is easy to see that the $n$th power of this permutation gives the identical transformation, so we
obtain the last {\em braid group relation}.

\begin{lemma} \label{lem-braid-2}
For any $n\ge3$, we have the second {\em braid group relation}
\begin{gather*}
\bigl(R_{n-1,n}R_{n-2,n-1}\cdots R_{2,3}R_{1,2}\bigr)^n=\hbox{\rm Id}.
\end{gather*}
\end{lemma}

\subsubsection[Braid group relations for geodesic functions of $D_n$-algebras]{Braid group relations for geodesic functions of $\boldsymbol{D_n}$-algebras}

It is possible to express readily the action of the braid group on
the level of the geodesic functions  $G_{i,j}$, $i,j=1,\dots,n$,
interpreted also as entries of the $n\times n$-matrix ${\mathcal D}$
(the elements that are not indicated remain invariant):
\begin{gather}
R_{i,i+1}{\mathcal D}={\tilde{\mathcal D}},\  \ \ \hbox{where} \
\left\{
\begin{array}{ll}
  {\tilde G}_{i+1,k}=G_{i,k}, & k\ne i,i+1,\\
  {\tilde G}_{i,k}=G_{i,k}G_{i,i+1}-G_{i+1,k}, &  k\ne i,i+1,\\
  {\tilde G}_{k,i+1}=G_{k,i}, &  k\ne i,i+1, \\
  {\tilde G}_{k,i}=G_{k,i}G_{i,i+1}-G_{k,i+1}, &  k\ne i,i+1, \\
  {\tilde G}_{i,i+1}=G_{i,i+1}, &  \\
  {\tilde G}_{i+1,i+1}=G_{i,i}, &  \\
  {\tilde G}_{i,i}=G_{i,i}G_{i,i+1}-G_{i+1,i+1}, &  \\
  {\tilde G}_{i+1,i}=G_{i+1,i}+G_{i,i+1}G_{i,i}^2-2G_{i,i}G_{i+1,i+1}. &  \\
\end{array}%
\right.\hspace{-5mm}
\label{Dn-braid-cl}
\end{gather}
The f\/irst braid group relation follows in this case as well from the three-step process,
but it can be verif\/ied explicitly that the following lemma holds.

\begin{lemma} \label{lem-braid-Dn}
For any $n\ge3$, we have the {\em braid group relation} for transformations \eqref{Dn-braid-cl}:
\begin{gather*}
R_{i-1,i}R_{i,i+1}R_{i,i-1}{\mathcal D}=R_{i,i+1}R_{i-1,i}R_{i,i+1}{\mathcal D},\qquad 2\le i\le n-1.
\end{gather*}
\end{lemma}

Note that the second braid-group relation (see
Lemma~\ref{lem-braid-2}) is lost in the case of $D_n$-algebras.

To present the braid-group action in the matrix-action (covariant)
form (\ref{BAB}) note that the combinations $G_{k,j}$, $G_{j,k}$,
and $G_{k,k}G_{j,j}$ have similar transformation laws in
(\ref{Dn-braid-cl}) in the case where at least one of the indices
$j$ and $k$ is neither $i$ nor $i+1$, so we can try to construct
globally covariantly transformed matrices from linear combinations
of the above (coef\/f\/icients of these combinations can be dif\/ferent
above and below the diagonal). Note that (since the braid-group
transformation acts on the $A_n$ subgroup of $D_n$ in the same way
as before), the matrices ${\mathcal A}$ (\ref{A-matrix}) and
${\mathcal A}^T$ are transformed as in (\ref{BAB}); the analysis
shows that we also have two new matrices, ${\mathcal R}$ and
$\mathcal S$, with the {\em same} transformation law as in
(\ref{BAB}):
\begin{gather}
\label{R.cl}
({\mathcal R})_{i,j}=\left\{
\begin{array}{cc}
                        G_{j,i}+G_{i,j}-G_{i,i}G_{j,j} &\quad j>i, \\
                        -G_{j,i}-G_{i,j}+G_{i,i}G_{j,j} &\quad  j<i, \\
                               0 &\quad  j=i,
                             \end{array}
                             \right.
\\
\label{S.cl}
({\mathcal S})_{i,j}=G_{i,i}G_{j,j}\qquad \hbox{for all}\quad 1\le i,j\le n,
\end{gather}
where ${\mathcal R}$ is skewsymmetric
(${\mathcal R}^T=-{\mathcal R}$) and ${\mathcal S}$ is symmetric
(${\mathcal S}^T={\mathcal S}$).

\begin{lemma} \label{lem-braid-Dn-matrix}
Any linear combination $w_1{\mathcal A}+w_2{\mathcal A}^T+\rho
{\mathcal R}+\sigma {\mathcal S}$ with complex $w_1$, $w_2$, $\rho$,
and $\sigma$ transforms in accordance with formula \eqref{BAB} under the braid-group
action.
\end{lemma}

We postpone the discussion of modular invariants constructed from
these four matrices till the discussion of the quantum $D_n$
braid-group action in Section~\ref{sss:q-Dn}.

\section{Quantum Teichm\"uller spaces of windowed surfaces}\label{ss:q}

\subsection{Canonical quantization of the Poisson algebra}

A quantization of a Poisson manifold, which is equivariant under the action of a discrete
group~$\cal D$,
is a family of $*$-algebras ${\cal A}^\hbar$ depending on a positive
real parameter $\hbar$ with
$\cal D$ acting by outer automorphisms and having the
following properties:

\begin{itemize}\itemsep=0pt

\item [{\bf 1.}] (Flatness.)  All algebras are isomorphic (noncanonically)
as linear spaces.

\item [{\bf 2.}]
(Correspondence.) For $\hbar=0$, the algebra is isomorphic as a $\cal
D$-module to the $*$-algebra of complex-valued functions on the Poisson
manifold.

\item [{\bf 3.}]  (Classical Limit.) The Poisson bracket on ${\cal A}^0$ given
by $\{a_1, a_2\} = \lim\limits_{\hbar \rightarrow 0}\frac{[a_1,a_2]}{\hbar}$ coincides with the
Poisson bracket given by the Poisson structure of the manifold.

\end{itemize}

Fix a cubic fatgraph $\Gamma_{g,\delta}$ as a spine of $\Sigma_{g,\delta}$,
and let ${\cal T}^\hbar={\cal T}^\hbar(\Gamma_{g,\delta})$
be the algebra generated by $Z_\alpha^\hbar$, one generator for each
unoriented edge $\alpha$ of $\Gamma_{g,\delta}$, with relations
\begin{gather}
\label{comm}
[Z^\hbar_\alpha, Z^\hbar_\beta ] = 2\pi i\hbar\{Z_\alpha, Z_\beta\}
\end{gather}
(cf.\ (\ref{WP-PB})) and the
$*$-structure
\[
(Z^\hbar_\alpha)^*=Z^\hbar_\alpha,
\]
where $Z_\alpha$  and
$\{\cdot,\cdot\}$ denotes the respective coordinate functions
and the Poisson bracket on the classical Teichm\"uller space.
Because of (\ref{WP-PB}), the right-hand side of (\ref{comm}) is a constant
taking only f\/ive values $0$, \ $\pm 2\pi i \hbar$, and $\pm 4 \pi i \hbar$
depending upon f\/ive variants of identif\/ications of endpoints of edges labelled $\alpha$ and $\beta$.

All the standard statements that we have in the case of Teichm\"uller spaces of
nonwindowed Riemann surfaces are transferred to the windowed surface case.

\begin{lemma}\label{nondeg}
The center ${\cal Z}^\hbar$ of the algebra ${\cal T}^\hbar$ is generated by the sums
$\sum\limits_{\alpha\in I}{Z^\hbar_\alpha}$ over all edges $\alpha\in I$ surrounding
a given boundary component, the center has dimension $s$, and the Poisson structure is
nondegenerate on the quotient ${\cal T}^\hbar/{\cal Z}^\hbar$.
\end{lemma}

The examples of such {\em boundary-parallel} curves are again in Fig.~\ref{fi:center}. Of course,
those are the same curves that provide the center of the Poisson algebra.

A standard Darboux-type theorem for nondegenerate Poisson structures then gives
the following result.

\begin{corollary}\label{corr11}
There is a basis for ${\cal T}^\hbar/{\cal Z}^\hbar$
given by operators $p_i$, $q_i$, for $i=1,\ldots ,3g-3+s+\sum\limits_{j=1}^s\delta_j$
satisfying the standard commutation relations $[p_i,q_j]=2\pi i \hbar\delta_{ij}$.
\end{corollary}

Now, def\/ine the Hilbert space ${\cal H}$ to be the set of all $L^2$ functions in the
$q$-variables and let each $q$-variable act by multiplication and each corresponding
$p$-variable act by dif\/ferentiation, $p_i=2\pi i\hbar\,{{\partial}\over{\partial q_i}}$. For
dif\/ferent choices of diagonalization of non-degenerate Poisson structures, these
Hilbert spaces are canonically isomorphic.

\subsection[Quantum flip transformations]{Quantum f\/lip transformations}

The Whitehead move becomes now a morphism of (quantum) algebras.
The {\it quantum Whitehead move} or {\it flip} along an edge of $\Gamma$ by equation (\ref{abc}) is
described by the (quantum) function~\cite{ChF}
\begin{gather} \label{phi}
\phi(z)\equiv \phi^\hbar(z) =
-\frac{\pi\hbar}{2}\int_{\Omega} \frac{e^{-ipz}}{\sinh(\pi p)\sinh(\pi \hbar
p)}dp,
\end{gather}
where the contour $\Omega$ goes along
the real axis bypassing the origin from above.
For each unbounded self-adjoint operator $Z^\hbar$ on ${\cal H}$, $\phi^\hbar (Z^\hbar)$ is a well-def\/ined
unbounded self-adjoint operator on ${\cal H}$.

The function $\phi^\hbar(Z)$ satisf\/ies the relations (see~\cite{ChF})
\begin{gather*}
\phi^\hbar(Z)-\phi^\hbar(-Z)=Z,
\nonumber
\\
\phi^\hbar(Z+i\pi\hbar)-\phi^\hbar(Z-i\pi\hbar)=\frac{2\pi i\hbar}{1+\e^{-Z}},
\nonumber
\\
\phi^\hbar(Z+i\pi)-\phi^\hbar(Z-i\pi)=\frac{2\pi i}{1+\e^{-Z/\hbar}}
\nonumber
\end{gather*}
and is meromorphic in the complex plane with the poles at the
points $\{\pi i(m+n\hbar),\ m,n\in {\Bbb Z}_+\}$ and
$\{-\pi i(m+n\hbar),\ m,n\in {\Bbb Z}_+\}$.

The function $\phi^\hbar(Z)$ is therefore holomorphic in the strip
$|\hbox{Im\,}Z|<\pi\hbox{\,min\,}(1,\hbox{Re\,}\hbar)-\epsilon$ for any $\epsilon>0$,
so we need only its asymptotic behavior as
$Z\in{\Bbb R}$ and $|Z|\to\infty$, for which we have (see, e.g., \cite{Kashaev3})
\begin{gather*}
\phi^\hbar(Z)\bigl|_{|Z|\to\infty}=(Z+|Z|)/2+O(1/|Z|).
\end{gather*}

We then have the following theorem \cite{ChF,Kashaev}

\begin{theorem}\label{th-Q}
The family of algebras ${\cal T}^\hbar={\cal T}^\hbar(\Gamma_{g,\delta} )$ is a
quantization of ${\cal T}^H_{g,\delta}$ for any cubic fatgraph spine
$\Gamma_{g,\delta}$ of $\Sigma_{g,\delta}$, that is,
\begin{itemize}\itemsep=0pt
\item In the limit $\hbar \mapsto 0$, morphism \eqref{abc} using
\eqref{phi} coincides with classical morphism \eqref{abc} with $\phi(Z)=\log(1+\e^Z)$.
\item Morphism \eqref{abc} using \eqref{phi} is indeed a morphism of $*$-algebras.
\item A flip $W_Z$ satisfies $W_Z^2=I$,
\eqref{abc}, and flips satisfy the commutativity relation.
\item Flips satisfy the pentagon relation.
\item The morphisms ${\cal T}^\hbar(\Gamma)\rightarrow {\cal
T}^{1/\hbar}(\Gamma)$ given by $Z^\hbar_\alpha \mapsto Z^{1/\hbar}_\alpha$
commute with morphisms \eqref{abc}.
\end{itemize}
\end{theorem}

\subsection{Geodesic length operators}

We next embed the algebra of geodesic functions
(\ref{G}) into a suitable
completion of the constructed algebra ${\cal T}^\hbar$.
For any $\gamma$,
the geodesic function $G_\gamma$ can be expressed  in terms of
shear coordinates on ${\cal T}^H$:
\begin{gather}
G_\gamma \equiv \tr P_{Z_1\cdots Z_n}=
\sum_{j\in J}\exp\left\{{\frac{1}{2} \sum_{\alpha \in E(\Gamma)}
m_j(\gamma,\alpha) Z_\alpha}\right\},  \label{clen}
\end{gather}
where $m_j(\gamma,\alpha)$ are
integers and $J$ is a f\/inite set of indices.

In general, sets of integers $\left\{m_j(\gamma,\alpha)\right\}_{\alpha=1}^{6g-6+3s+2\#\delta}$
may coincide for dif\/ferent $j_1,j_2\in J$; we however distinguish between them as soon
as they come from dif\/ferent products of exponentials $e^{\pm Z_i/2}$ in
traces of matrix products in (\ref{clen}).

For any closed path $\gamma$ on $\Sigma_{g,\delta}$, def\/ine the {\it quantum geodesic}
operator $G^\hbar_\gamma \in {\cal T}^\hbar$ to be
\begin{gather} \label{qlen}
G^\hbar_\gamma \equiv
\ORD{\tr P_{Z_1\dots Z_n}}\equiv
\sum_{j\in J}
\exp\left\{{\frac{1}{2} \sum_{\alpha \in E(\Gamma_{g,\delta})}
\bigl(m_j(\gamma,\alpha) Z^\hbar_\alpha
+2\pi i\hbar
c_j(\gamma,\alpha)
\bigr)}\right\},
\end{gather}
where the {\it quantum ordering} $\ORD{\cdot}$ implies that we vary the classical
expression (\ref{clen}) by introducing additional
integer coef\/f\/icients $c_j(\gamma,\alpha)$,
which must be determined from the conditions below.

That is, we assume that each term in the classical expression (\ref{clen}) can get
multiplicative corrections only of the form $q^n$, $n\in {\Bbb Z}$, with
\begin{gather*}
q\equiv e^{-i\pi\hbar}.
\end{gather*}

We often call a quantum geodesic function merely a
quantum geodesic because quantum objects admit only a functional description.

We now formulate the def\/ining properties of quantum geodesics.
\begin{itemize}\itemsep=0pt
\item[{\bf 1.}] If closed paths $\gamma$ and $\gamma^\prime$ do not
intersect, then the operators $G^\hbar_\gamma$ and
$G^\hbar_{\gamma^\prime}$ commute.
\item[{\bf 2.}]  {\it Naturality.} The mapping class group
$MC(\Sigma_{g,\delta})$ (\ref{abc}) acts naturally, i.e., for any $\{G^\hbar_\gamma\}$,
$W^\hbar \in MC(\Sigma_{g,\delta})$, and closed path $\gamma$ in a spine $\Gamma_{g,\delta}$ of
$\Sigma_{g,\delta}$, we have
\[
W^\hbar(G^\hbar_\gamma) = G^\hbar_{W(\gamma)}.
\]
\item[{\bf 3.}]  {\it Geodesic algebra}. The product of two quantum geodesics
is a linear combination of quantum multicurves governed by the (quantum) skein
relation below.
\item[{\bf 4.}] {\it Orientation invariance.} Quantum traces of
direct and inverse geodesic operators coincide.
\item[{\bf 5.}] {\it Exponents of geodesics.}
A quantum geodesic $G^\hbar_{n\gamma}$ corresponding to the $n$-fold concatenation of $\gamma$
is expressed via $G^\hbar_{\gamma}$ exactly as in the classical case, namely,
\begin{gather}
\label{cheb}
G^\hbar_{n\gamma}=2T_n\bigl(G^\hbar_{\gamma}/2\bigr),
\end{gather}
where $T_n(x)$ are Chebyshev's polynomials.
\item[{\bf 6.}] {\it Hermiticity.} A quantum geodesic is a Hermitian
operator having by def\/inition a real spectrum.
\end{itemize}

We shall let  the standard normal ordering symbol ${:}\e^{a_1}\e^{a_2}\cdots\e^{a_n}{:}$ denote the {\it Weyl
ordering} $\e^{a_1+\cdots+a_n}$,  i.e.,
\begin{gather*}
{:}\e^{a_1}\e^{a_2}\cdots\e^{a_n}{:}=1+(a_1+\cdots +a_n)+{1\over {2!}}~(a_1+\cdots +a_n)(a_1+\cdots
+a_n)+\cdots
\end{gather*}
for any set of exponents with $a_i\neq -a_j$ for $i\neq j$,
In particular, the Weyl ordering implies total symmetrization in the subscripts.

We have~\cite{ChF} the proposition, which can be extended to the case of windowed surfaces
assuming the modif\/ication of the ``old" notion of graph simple geodesics.

\begin{definition}
For a spine $\Gamma_{g,\delta}$, we call a geodesic {\it graph simple} if it does not pass twice
through any of inner edges of the graph and has at most one inversion at any of pending edges.
\end{definition}

\begin{proposition}\label{lem31}
For any graph simple geodesic $\gamma$ with respect to any spine $\Gamma$,
the coefficients $c_j(\gamma,\alpha)$ in \eqref{qlen} are identically zero, i.e.,
the quantum ordering is the Weyl ordering.
\end{proposition}

\begin{proof} Let us again denote by $Y^\hbar_i$, $i=1,\dots, 6g-6+3s+\#\delta$,
the quantum shear coordinates of inner edges and by $Z^\hbar_j$, $j=1,\dots, \#\delta$
the quantum shear coordinates of pending edges. But the latter always come in the
combination $X_{Z^\hbar_j}F X_{Z^\hbar_j}=X_{2Z^\hbar_j}$, so, considering
term-by-term the trace of the matrix product for a quantum graph simple geodesic, we f\/ind that
we can expand it in Laurent monomials in $\e^{Y^\hbar_i/2}$ and $\e^{Z^\hbar_j}$.
It is easy to see that each term $\e^{Y^\hbar_i/2}$ and $\e^{Z^\hbar_j}$ comes either in power $+1$, or
$-1$ in the corresponding monomial and there are no equivalent monomials in
the sum.  This means that in order to have a {Hermitian} operator, we must apply the
Weyl ordering with no additional $q$-factors (by the correspondence principle, each
such factor must be again a Laurent monomial in $q$ standing by the corresponding
term, which breaks the self-adjointness unless all such monomials are unity).  Since quantum
Whitehead moves must preserve the property of being Hermitian, if a graph-simple geodesic transforms
to another graph-simple geodesic, then a Weyl-ordered expression transforms to a Weyl-ordered
expression, and only these expressions are self-adjoint.
\end{proof}

\begin{example}
For the $A_3$ algebra graph in Fig.~\ref{fi:An}, we have exactly three graph simple
geodesics with the corresponding geodesic functions $G^\hbar_{12}$, $G^\hbar_{23}$, and $G^\hbar_{13}$ given
by formulas (\ref{G12}) (which are written already in the Weyl-ordered form), and if
we consider, for instance, the product
\begin{gather}
\label{TurPR}
G^\hbar_{23}G^\hbar_{12}=q^{-1}G^\hbar_{1232}+q G^\hbar_{13},
\end{gather}
where $G^\hbar_{1232}$ and $G^\hbar_{13}$ correspond to respective cases (a) and (b) of resolving crossing
of the geodesics $\gamma_{23}$ and $\gamma_{12}$ near the dot-vertex~2 in Fig.~\ref{fi:dot-skein}.
Note that $G^\hbar_{13}$ is also Weyl-ordered, $G_{13}=e^{Z_3+Z_1}+e^{Z_3-Z_1}+e^{-Z_3-Z_1}$
whereas
\begin{gather*}
G^\hbar_{1232}=e^{Z_1+2Z_2+Z_3}+e^{Z_1+2Z_2-Z_3}+e^{Z_1-2Z_2-Z_3}+e^{-Z_1-2Z_2-Z_3}
+(q^2+q^{-2})e^{Z_1-Z_3}
\end{gather*}
apparently is not Weyl-ordered.
Product of the same operators as in (\ref{TurPR}) but taken in opposite order gives
$G^\hbar_{12}G^\hbar_{23}=qG^\hbar_{1232}+q^{-1} G^\hbar_{13}$, so, introducing the $q^2$-commutator
$[A,B]_{q^2}\equiv qAB-q^{-1}BA$ and $\xi=q^2-q^{-2}$, we have the quantum $A_3$-algebra:
\begin{gather*}
[G^\hbar_{23},G^\hbar_{12}]_{q^2}=\xi G^\hbar_{13},\qquad
[G^\hbar_{13},G^\hbar_{23}]_{q^2}=\xi G^\hbar_{12},\qquad
[G^\hbar_{12},G^\hbar_{13}]_{q^2}=\xi G^\hbar_{23}.
\end{gather*}
\end{example}

\subsection{Quantum skein relations}\label{ss:skein}

We now formulate the {\em general} rules that allow one to disentangle the product of
any two quantum geodesics.

Let $G^\hbar_1$ and $G^\hbar_2$ be two quantum geodesic operators corresponding to
geodesics $\gamma_1$ and $\gamma_2$ where all the inversion relations are resolved using
the dot-vertex construction (see Fig.~\ref{fi:dot-vertex}). Then
\begin{itemize}\itemsep=0pt
\item
We must apply the {\em quantum skein relation}\footnote{Here the order of crossing lines corresponding to
$G^\hbar_1$ and $G^\hbar_2$ depends on which quantum geodesic
occupies the f\/irst place in the product;
the rest of the graph remains unchanged for all items in (\ref{skein}).}
\begin{gather}
\begin{gathered}[c]
\includegraphics{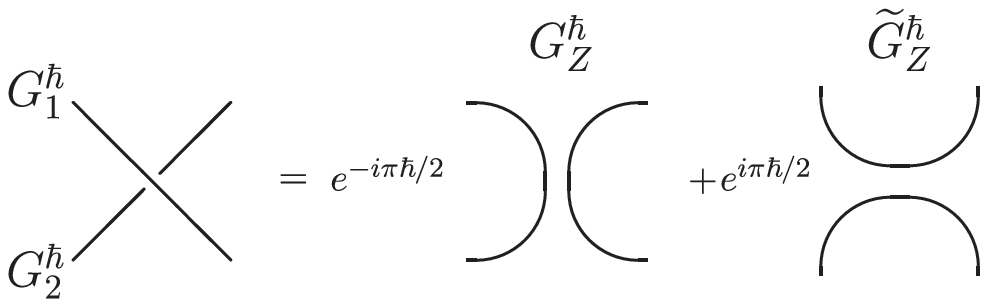}
\end{gathered}
\label{skein}
\end{gather}
{\em simultaneously} at {\em all} intersection points.
\item After the application of the quantum skein relation we can
obtain empty (contractible) loops; we assign the factor $-q-q^{-1}$
to each such loop and this suf\/f\/ices to ensure the quantum Riedemeister moves.
\item We can also obtain loops that are homeomorphic to going around a dot-vertex;
as in the classical case, we claim the corresponding geodesic functions to vanish,
so we disregard all such cases of geodesic laminations in the quantum case as well.
\end{itemize}

The main lemma is in order.

\begin{lemma}[\cite{ChF,Teshner}]\label{lem34}
There exists
a unique quantum ordering $\ORD{\cdots}$ \eqref{qlen}, which is
generated by the quantum geodesic
algebra \eqref{skein} and is consistent with the quantum mapping
class groupoid transformations \eqref{abc}, i.e., so that the quantum
geodesic algebra is invariant under the action of the quantum mapping
class groupoid.
\end{lemma}

\subsection{Quantum braid group relation}\label{sss:q-braid}

\subsubsection[Quantum $A_n$-algebra]{Quantum $\boldsymbol{A_n}$-algebra}\label{sss:q-braid-An}

We now consider the quantum geodesic functions associated with paths in the $A_n$-algebra
pattern in Fig.~\ref{fi:An}.
From the quantum skein relation, it is easy to obtain quantum transformations for the quantum
geodesic functions $G^\hbar_{i,j}$. We introduce the ${\mathcal A}^\hbar$-matrix
\begin{gather}
\label{A-matrix-q}
{\mathcal A}^\hbar=\left(\begin{array}{ccccc}
                     q & G^\hbar_{1,2} & G^\hbar_{1,3} & \dots & G^\hbar_{1,n} \\
                     0 & q & G^\hbar_{2,3} & \dots & G^\hbar_{2,n} \\
                     0 & 0 & q & \ddots & \vdots \\
                     \vdots & \vdots & \ddots & \ddots & G^\hbar_{n-1,n} \\
                     0 & 0 & \dots & 0 & q \\
                   \end{array}
\right)
\end{gather}
associating the Hermitian operators $G^\hbar_{i,j}$ with the quantum
geodesic functions. Using the skein relation,
we can then present the action of the braid group element $R^\hbar_{i,i+1}$ exclusively in terms
of the geodesic functions from this, f\/ixed set:
$R^\hbar_{i,i+1}{\mathcal A}^\hbar={\tilde{\mathcal A}^\hbar}$, where
\begin{gather}
\label{R-matrix-q}
\begin{array}{ll}
  {\tilde G}^\hbar_{i+1,j}=G^\hbar_{i,j}, & j>i+1,\\
  {\tilde G}^\hbar_{j,i+1}=G^\hbar_{j,i}, & j<i, \\
  {\tilde G}^\hbar_{i,j}=qG^\hbar_{i,j}G^\hbar_{i,i+1}-q^{2}G^\hbar_{i+1,j}
  =q^{-1}G^\hbar_{i,i+1}G^\hbar_{i,j}-q^{-2}G^\hbar_{i+1,j}, & j>i+1, \\
  {\tilde G}^\hbar_{j,i}=qG^\hbar_{j,i}G^\hbar_{i,i+1}-q^2G^\hbar_{j,i+1}
  =q^{-1}G^\hbar_{i,i+1}G^\hbar_{j,i}-q^{-2}G^\hbar_{j,i+1}, & j<i, \\
  {\tilde G}^\hbar_{i,i+1}=G^\hbar_{i,i+1}. &  \\
\end{array}
\end{gather}
We can again present this transformation via the special matrices
$B^\hbar_{i,i+1}$ of the block-diagonal form
\begin{gather*}
B^\hbar_{i,i+1}=\begin{array}{c}
            \vdots \\
            i \\
            i+1 \\
            \vdots \\
          \end{array}
          \left(
          \begin{array}{cccccccc}
            1 &  &  &  &  &  &  &  \\
             & \ddots &  &  &  &  &  &  \\
             &  & 1 &  &  &  &  &  \\
             &  &  &   q^{-1}G^\hbar_{i,i+1} & -q^{-2} &  & &  \\
            &  &  & 1 & 0 &  &  &  \\
             &  &  &  &  & 1 &  &  \\
             &  &  &  &  &  & \ddots &  \\
             &  &  &  &  &  &  & 1 \\
          \end{array}
          \right).
\end{gather*}
Then, the action of the quantum braid group generator $R^\hbar_{i,i+1}$ on ${\mathcal A}^\hbar$
can be expressed as the matrix product (taking into account the noncommutativity of quantum matrix entries)
\begin{gather}
\label{BAB-q}
R^\hbar_{i,i+1}{\mathcal A}^\hbar=B^\hbar_{i,i+1}{\mathcal A}^\hbar\bigl(B^\hbar_{i,i+1}\bigr)^{\dagger}
\end{gather}
with $\bigl(B^{\hbar}_{i,i+1}\bigr)^{\dagger}$ the matrix Hermitian conjugate to
$B^\hbar_{i,i+1}$ (its nontrivial $2\times 2$-block has the form $\left(%
\begin{array}{cc}
  qG^\hbar_{i,i+1} & 1 \\
  q^2 & 0 \\
\end{array}%
\right)$).
Using the same technique as above, it is then straightforward to prove the following lemma.

\begin{lemma} \label{lem-q-braid}
For any $n\ge3$, we have the {\em quantum braid group relations}
\begin{gather}
\label{RRR-q}
R^\hbar_{i-1,i}R^\hbar_{i,i+1}R^\hbar_{i-1,i}=R^\hbar_{i,i+1}R^\hbar_{i-1,i}R^\hbar_{i,i+1},\qquad 2\le i\le n-1,
\\
\label{Rn-q}
\bigl(R^\hbar_{n-1,n}R^\hbar_{n-2,n-1}\cdots R^\hbar_{2,3}R^\hbar_{1,2}\bigr)^n=\hbox{\rm Id}.
\end{gather}
\end{lemma}

\subsubsection[Quantum $D_n$-algebra]{Quantum $\boldsymbol{D_n}$-algebra}\label{sss:q-Dn}

\begin{figure}[tb]
\centerline{\includegraphics{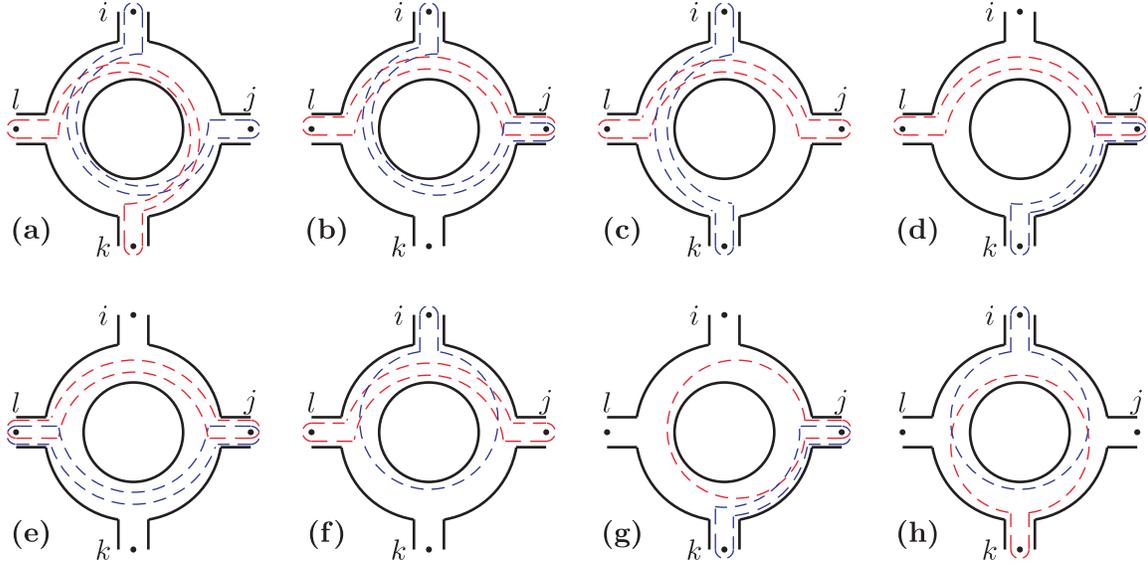}}
\caption{Eight cases of nontrivial intersections of geodesics from
the set $G_{ij}$, $i,j=1,\dots,n$ in the case of the $D_n$-algebra.}
\label{fi:Dn-variants}
\end{figure}

We now quantize the Poisson algebra of geodesic functions $G_{ij}$ corresponding to paths
as shown in Fig.~\ref{fi:Dn}. We have there eight possible variants of nontrivial intersections
shown in Fig.~\ref{fi:Dn-variants}.

The corresponding quantum permutation relations read\footnote{Deriving these relations requires a tedious combinatorial
analysis based on quantum skein relations formulated in Section~\ref{ss:skein}.}
($q=e^{-i\pi\hbar}$, $\xi\equiv q^2-q^{-2}$)
\begin{gather}
\hbox{Case (a)}\qquad [G^\hbar_{ij},G^\hbar_{kl}]=\xi\biggl(G^\hbar_{kj}G^\hbar_{li}-G^\hbar_{jk}G^\hbar_{il}
+G^\hbar_{jl}G^\hbar_{ik}-G^\hbar_{lj}G^\hbar_{ki}
\nonumber
\\
\phantom{\hbox{Case (a)}\qquad [G^\hbar_{ij},G^\hbar_{kl}]=}{}
+(q+q^{-1})(G^\hbar_{il}G^\hbar_{jj}G^\hbar_{kk}-G^\hbar_{kj}G^\hbar_{ll}G^\hbar_{ii})
\biggr);
\nonumber
\\
\hbox{Case (b)}\qquad qG^\hbar_{jl}G^\hbar_{ij}-q^{-1}G^\hbar_{ij}G^\hbar_{jl}=\xi\biggl(2G^\hbar_{il}-G^\hbar_{li}
-G^\hbar_{il}(G^\hbar_{jj})^2\biggr)
\nonumber
\\
\phantom{\hbox{Case (b)}\qquad qG^\hbar_{jl}G^\hbar_{ij}-q^{-1}G^\hbar_{ij}G^\hbar_{jl}=}{}
+\xi(q+q^{-1})G^\hbar_{ii}G^\hbar_{ll} +(q-q^{-1})G^\hbar_{lj}G^\hbar_{ji};
\nonumber
\\
\hbox{Case (c)}\qquad [G^\hbar_{ik},G^\hbar_{jl}]=\xi\biggl(G^\hbar_{jk}G^\hbar_{il}-G^\hbar_{ji}G^\hbar_{lk}\biggr);
\nonumber
\\
\hbox{Case (d)}\qquad qG^\hbar_{jl}G^\hbar_{kj}-q^{-1}G^\hbar_{kj}G^\hbar_{jl}=\xi G^\hbar_{kl};
\label{Dn-q}
\\
\hbox{Case (e)}\qquad [G^\hbar_{jl},G^\hbar_{lj}]=\xi\biggl((G^\hbar_{ll})^2-(G^\hbar_{jj})^2\biggr);
\nonumber
\\
\hbox{Case (f)}\qquad [G^\hbar_{jl},G^\hbar_{ii}]=\xi\biggl(G^\hbar_{ji}G^\hbar_{ll}-G^\hbar_{il}G^\hbar_{jj}\biggr);
\nonumber
\\
\hbox{Case (g)}\qquad qG^\hbar_{jj}G^\hbar_{kj}-q^{-1}G^\hbar_{kj}G^\hbar_{jj}=\xi G^\hbar_{kk},\qquad
qG^\hbar_{jk}G^\hbar_{jj}-q^{-1}G^\hbar_{jj}G^\hbar_{jk}=\xi G^\hbar_{kk};
\nonumber
\\
\hbox{Case (h)}\qquad [G^\hbar_{ii},G^\hbar_{kk}]=(q-q^{-1})\bigl(G^\hbar_{ik}-G^\hbar_{ki}\bigr).
\nonumber
\end{gather}

Although these relations not only contain triple terms in the r.h.s.\ but also noncommuting terms (this is the
price for closing the algebra), they nevertheless establish the lexicographic ordering on the corresponding set
of quantum variables $\{G^\hbar_{ij}\}$.

\begin{lemma}\label{lem-Dn-quantum}
Permutation relations treated as an abstract algebra
postulated by \eqref{Dn-q} satisfy the commutation Jacobi identities.
\end{lemma}

The proof is tedious but straightforward calculations. Note that algebra (\ref{Dn-q}) is self-consistent even without relation
to geometry of modular spaces; the similar phenomenon was already observed in the case of $A_n$-algebras. As regarding the
classif\/ication of cluster algebras in~\cite{FST}, we produced the corresponding algebras for disc and annulus with arbitrary
number of marked points (in our approach, a punctured disc is just an annulus with one hole of zero perimeter; using the
isomorphism in Theorem~\ref{th-mark} we can move all marked points to one boundary component). We do not know however
as yet an example
of such an algebraically closed set for a disc with two punctures (holes); this case deserves a separate investigation.

We now provide the quantum version of the braid group transformations (\ref{Dn-braid-cl}).
They are similar to (\ref{R-matrix-q}), the only actual distinction is in the
corresponding $2\times 2$-block of the matrix~${\mathcal D}$. It is however that distinction
that prevents us from writing these transformations in a form similar to (\ref{BAB-q}).

We therefore have for the quantum braid group transformation for the $D_n$-algebra
\begin{gather}
\label{R-matrix-Dn-q}
\begin{array}{ll}
  {\tilde G}^\hbar_{i+1,k}=G^\hbar_{i,k}, & k\ne i,i+1,\\
  {\tilde G}^\hbar_{i,k}=qG^\hbar_{i,i+1}G^\hbar_{i,k}-q^{2}G^\hbar_{i+1,k}
  =q^{-1}G^\hbar_{i,k}G^\hbar_{i,i+1}-q^{-2}G^\hbar_{i+1,k}, & k\ne i,i+1, \\
  {\tilde G}^\hbar_{k,i+1}=G^\hbar_{k,i}, & k\ne i,i+1,\\
  {\tilde G}^\hbar_{k,i}=qG^\hbar_{i,i+1}G^\hbar_{k,i}-q^{2}G^\hbar_{k,i+1}
  =q^{-1}G^\hbar_{k,i}G^\hbar_{i,i+1}-q^{-2}G^\hbar_{k,i+1}, & k\ne i,i+1, \\
  {\tilde G}^\hbar_{i,i+1}=G^\hbar_{i,i+1}, \qquad
  {\tilde G}^\hbar_{i+1,i+1}=G^\hbar_{i,i}, &  \\
  {\tilde G}^\hbar_{i,i}=qG^\hbar_{i,i+1}G^\hbar_{i,i}-q^{2}G^\hbar_{i+1,i+1}
  =q^{-1}G^\hbar_{i,i}G^\hbar_{i,i+1}-q^{-2}G^\hbar_{i+1,i+1}, &  \\
  {\tilde G}^\hbar_{i+1,i}=G^\hbar_{i+1,i}+G^\hbar_{i,i}G^\hbar_{i,i+1}G^\hbar_{i,i}
  -q^{-1}G^\hbar_{i+1,i+1}G^\hbar_{i,i}-qG^\hbar_{i,i}G^\hbar_{i+1,i+1}. &  \\
\end{array}%
\end{gather}

\begin{lemma} \label{lem-q-braid-Dn}
For any $n\ge3$, we have the {\em quantum braid group relations} \eqref{RRR-q}
for transformations \eqref{R-matrix-Dn-q} of quantum operators subject to
quantum algebra \eqref{Dn-q}.
\end{lemma}

Again, it the second identity (\ref{Rn-q}) is lost in the case of $D_n$ algebras.

\subsubsection[Matrix representation for $D_n$-algebra and invariants]{Matrix representation for $\boldsymbol{D_n}$-algebra and invariants}
\label{sss:q-Dn-mat}

We now construct the quantum analogues of (\ref{R.cl}) and
(\ref{S.cl}).

\begin{lemma} \label{lem-q-Dn-matrix}
The following four matrices (with operatorial entries), together
with all their linear combinations, transform in accordance with the
quantum braid-group action \eqref{BAB-q}: ${\mathcal A}^\hbar$
\eqref{A-matrix-q}, $\bigl({\mathcal A}^\hbar\bigr)^\dagger$,
${\mathcal R}^\hbar$, and ${\mathcal S}^\hbar$, where
\begin{gather}
\label{R.q}
({\mathcal R}^\hbar)_{i,j}=\left\{
\begin{array}{cc}
                        G^\hbar_{j,i}+q^2G^\hbar_{i,j}-qG^\hbar_{i,i}G^\hbar_{j,j} &\quad j>i, \\
                        -G^\hbar_{i,j}-q^{-2}G^\hbar_{j,i}+q^{-1}G^\hbar_{i,i}G^\hbar_{j,j} &\quad  j<i, \\
                               0 &\quad  j=i,
                             \end{array}
                             \right.\qquad
\bigl({\mathcal R}^\hbar\bigr)^\dagger=-{\mathcal R}^\hbar,
\\
({\mathcal
S}^\hbar)_{i,j}=G^\hbar_{i,i}G^\hbar_{j,j}\qquad \hbox{for
all}\quad 1\le i,j\le n, \qquad \bigl({\mathcal S}^\hbar\bigr)^\dagger
={\mathcal S}^\hbar.
\nonumber
\end{gather}
\end{lemma}

\begin{remark}
We took particular form\footnote{Recall that we can ``play''
with coef\/f\/icients adding matrices ${\mathcal A}^\hbar$ and
$\bigl({\mathcal A}^\hbar\bigr)^\dagger$.} (\ref{R.q}) of the matrix
${\mathcal R}^\hbar$ because, in the case $n=2$, the combination
\begin{gather*}
G^\hbar_{1,1}G^\hbar_{2,2}-qG^\hbar_{1,2}-q^{-1}G^\hbar_{2,1}=
G^\hbar_{2,2}G^\hbar_{1,1}-q^{-1}G^\hbar_{1,2}-qG^\hbar_{2,1}
\end{gather*}
is a central element of the (quantum) algebra $D_2$; the other
central element is
\begin{gather*}
G^\hbar_{1,2}G^\hbar_{2,1}-q^2(G^\hbar_{2,2})^2-q^{-2}(G^\hbar_{1,1})^2=
G^\hbar_{2,1}G^\hbar_{1,2}-q^{-2}(G^\hbar_{2,2})^2-q^2(G^\hbar_{1,1})^2.
\end{gather*}
Also, exactly with such a combination, diagonal elements remain
zeros acquiring no $q$-corrections. And, eventually, only this
matrix possesses the (quantum) cyclic symmetry below.
\end{remark}

A cyclic permutation of indices $P:\,i\mapsto i+1 \mod n;\ j\mapsto
j+1 \mod n$ destroys the structure of the matrix ${\mathcal
A}^\hbar$ and results in the following transformations for
${\mathcal R}^\hbar$ and ${\mathcal S}^\hbar$:
\begin{gather}
P:\ {\mathcal R}^\hbar\mapsto
\left(%
\begin{array}{cccc}
  0 & 1 &  &  \\
   & \ddots & \ddots &  \\
   &  & \ddots & 1 \\
  -q^{-2} &  &  & 0 \\
\end{array}%
\right) {\mathcal R}^\hbar
\left(%
\begin{array}{cccc}
  0 &  &  & -q^2 \\
  1 & \ddots & &  \\
   & \ddots & \ddots &  \\
   &  & 1 & 0 \\
\end{array}%
\right),
\nonumber
\\
P:\ {\mathcal S}^\hbar\mapsto
\left(%
\begin{array}{cccc}
  0 & 1 &  &  \\
   & \ddots & \ddots &  \\
   &  & \ddots & 1 \\
  1 &  &  & 0 \\
\end{array}%
\right) {\mathcal S}^\hbar
\left(%
\begin{array}{cccc}
  0 &  &  & 1 \\
  1 & \ddots & &  \\
   & \ddots & \ddots &  \\
   &  & 1 & 0 \\
\end{array}%
\right).
\nonumber
\end{gather}
These transformations together with (\ref{BAB-q})
generate a full modular group. This means that, at least in the
classical case, $\det {\mathcal R}$ is the mapping-class group
invariant and lies therefore in the center of the Poisson algebra.
Same is true for ${\mathcal S}$, but $\det {\mathcal S}\equiv 0$
whereas $\det {\mathcal R}$ is nonzero for even $n=2m$ (and vanishes
for odd $n$): denoting $Q_{i,j}:=({\mathcal R})_{i,j}$ for $i<j$, we
have that $\det {\mathcal R}=\hbox{Pf\,}({\mathcal R})^2$, where the
Pfaf\/f\/ian $\hbox{Pf\,}({\mathcal R})$ is given by
the Grassmann-variable integral
\begin{gather*}
\hbox{Pf\,}({\mathcal R})=\int\cdots\int d\theta_1\cdots d\theta_{2m}e^{\sum\limits_{1\le i<j\le
2m}\theta_i Q_{i,j}\theta_j}
\end{gather*}
and is described by all possible (signed) pairings in the set
$\theta_1\theta_2\cdots\theta_{2m-1}\theta_{2m}$ where
$\langle\theta_i\theta_j\rangle=Q_{i,j}$ for $i<j$.

For example, for $m=2$, we have
\[
I_4=Q_{1,2}Q_{3,4}+Q_{1,4}Q_{2,3}-Q_{1,3}Q_{2,4}
\]
(recall that $Q_{i,j}=G_{i,j}+G_{j,i}-G_{i,i}G_{j,j}$ in the
classical case). In the quantum case, these elements obviously get
$q$-corrections to be calculated explicitly since the notion of a
quantum determinant is ambiguous; we hope to return to this problem
elsewhere.

Also, this construction provides just one central element of the
algebra $D_{2m}$; f\/inding other central elements in a regular way
(similar to the one in \cite{Bondal} where all the central elements
of the $A_n$-algebra were generated by $\det(\lambda^{-1}{\mathcal
A}+\lambda {\mathcal A}^T)$) needs further investigation.


\section{Multicurves and the double of windowed surface}\label{s:double}

\subsection{Multicurves for bordered surfaces}\label{ss:multi}

It is a standard notation that a multicurve (lamination)
is a collection of non(self)inter\-secting curves.
Apparently, in a new formulation with dot-vertices, this def\/inition can be literally
transferred to the case of surfaces with marked points on boundary components.

\begin{figure}[tb]
\centerline{\includegraphics{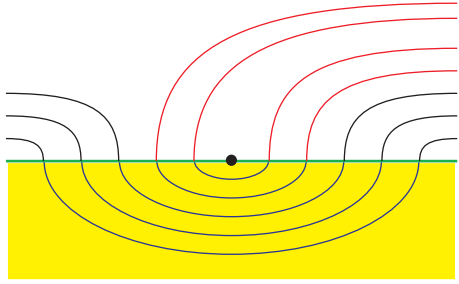}}
\caption{Closing the multicurve on a boundary component (shown by blue curves).
We let the bullet denote the dot-vertex. Lower (yellow) part is f\/ictitious and is
not present in the original geometry of surfaces with boundaries.}
\label{fi:bounding}
\end{figure}

\begin{definition} \label{def1}
Consider the homotopy class of a f\/inite collection
$C_e=\{ \gamma _1,\ldots ,\gamma _n\}$
of disjointly embedded (unoriented) simple curves $\gamma _i$ in a
topological windowed surface $\Sigma_{g,\delta}$. These curves are either
closed or terminate at windows. We impose the only restriction that we have
an even (can be zero) number of endpoints of these curves at each window.
An {\it even-based generalized multicurve} (eGMC) ${\mathcal C}_e$ in $\Sigma_{g,\delta}$
is a multiset based on $C$; we then have $s_i\geq 1$ parallel copies of components
of $C$, or in other words, positive integral weights $s_i$ on each component of $C$,
where $s_i$ is the multiplicity of $\gamma _i$ in ${\mathcal C}_e$. Further, given a
hyperbolic structure on $\Sigma_{g,\delta}$, we associate to~${\mathcal C}_e$ the product
\begin{gather*}
G_{{\mathcal C}_e}=G_{\gamma_1}^{s_1}\cdots G_{\gamma_n}^{s_n}
\end{gather*}
of geodesic operators (\ref{G}) of
all geodesics constituting a GMC; these operators Poisson commute in the classical
case since the components of $C$ are disjoint.

We therefore extend the standard def\/inition by introducing {\em nonclosed} curves that terminate at the
boundary components. The {\em only} restriction we impose is the {\em evenness} of
the number of curves terminating at each connected part of a boundary component
(passing through an inversion
geodesic line). Now in order to obtain the multicurve,
we connect pairwise the ends of these curves ({\em below} the inversion, or bounding,
geodesic line) as shown in Fig.~\ref{fi:bounding} obtaining therefore the collection of
topologically {\em closed} curves.

Each curve that is homeomorphic to the path encircling a single dot-vertex (see, e.g., case (d) in
Fig.~\ref{fi:dot-skein}) has {\em zero} geodesic function $G=\tr F$, and we therefore exclude all
such multicurves. At the same time, as before, a eGMC containing a
contractible component (of length zero) is minus twice the eGMC with this curve removed.
\end{definition}

\subsection{The double of the Riemann surface}\label{ss:double}

In order for our procedure to have a geometrical sense,
we now construct the double of the Riemann surface and transfer to this double the
structure of multicurves. Note that it is often
easier to prove properties of the geodesic function using the structure of the doubled Riemann
surface without boundary components.

\begin{figure}[tb]
\centerline{\includegraphics{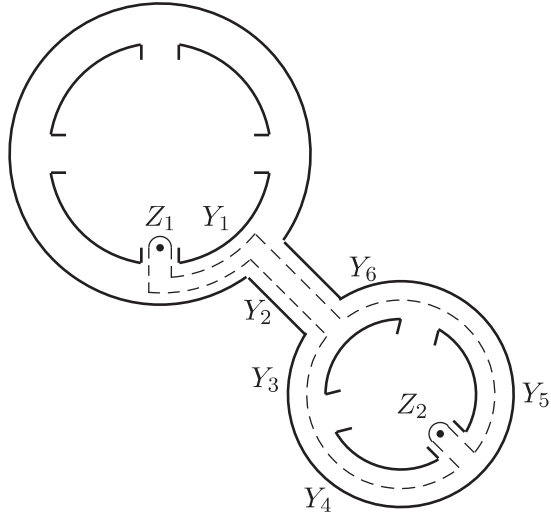}}
\caption{An example of geodesic curve with even number of inversions.}
\label{fi:geodesic}
\end{figure}

\begin{figure}[tb]
\centerline{\includegraphics{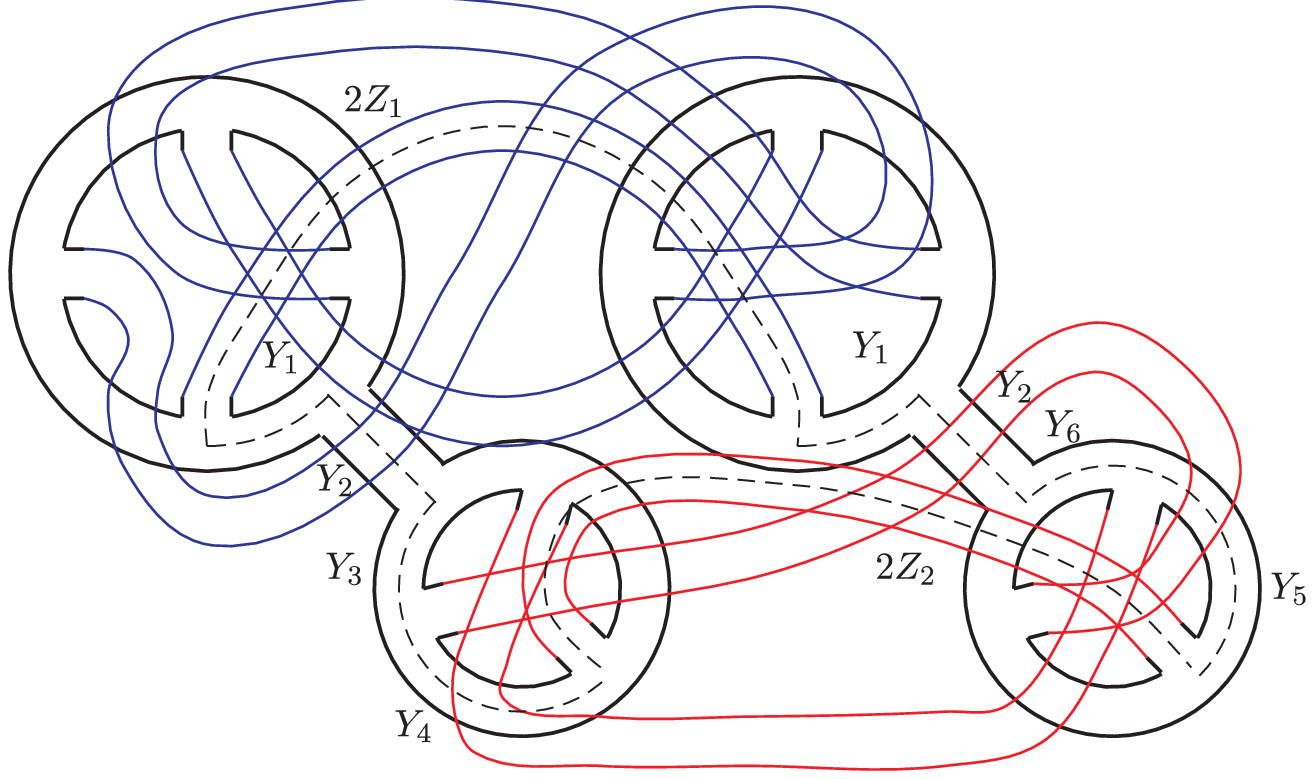}}
\caption{The double representation for the Riemann surface depicted in Fig.~\ref{fi:geodesic}.
Pending propagators on boundary components
are connected, each carrying the doubled shear coordinate $2Z_i$, in a crosswise manner.
We then obtain an oriented surface whose genus is given by (\ref{genus-double}). To each nonself\/intersecting
geodesic line that has an even number of inversions in the original picture, we may set into the correspondence
a unique (nonself\/intersecting) line on the doubled Riemann surface (dashed line for the geodesic line in
Fig.~\ref{fi:geodesic}).}
\label{fi:double}
\end{figure}

We f\/irst describe the doubled surface itself. Nonintersecting
curves must remain nonintersecting in the picture of double when we replace inversions merely
by (doubled) variables $2Z_i$ on the corresponding edges. To attain this, let us consider the example
of a geodesic line in Fig.~\ref{fi:geodesic}. If we clone this surface, drift apart
the copy of it and then connect all the pending ends with their twins preserving the orientation (as shown in Fig.~\ref{fi:double}),
then we obtain the new Riemann surface without pending edges (and, correspondingly, without windows).
Having the initial Riemann surface $\Sigma_{g,\delta}$ where $\delta$ is from (\ref{delta}),
$s$ the number of boundary components, and~$\delta_j$ the number of marked points on
the $j$th component (can be zero), we obtain that upon the joining as shown in Fig.~\ref{fi:double}, each boundary component
with {\em even} $\delta_j$ (including the zero one) generates two boundary components without marked points (holes)
in the doubled Riemann surface, whereas each boundary component with {\em odd} $\delta_j$ generates exactly one hole in
the new doubled Riemann surface. Then, easy calculation using the Euler characteristic formula gives the answer
for the genus $\hat g$ and number of holes $\hat s$ of the doubled Riemann surface $\Sigma_{\hat g,\hat s}$:
\begin{gather}
\label{genus-double}
\hat g=2g-1+\frac12 \left(\sum_{j=1}^s\delta_j+\#\{\hbox{odd $\delta_j$}\}\right),\qquad \hat s=2s-\#\{\hbox{odd $\delta_j$}\}.
\end{gather}

Note that while f\/lip morphisms on inner edges of the initial graph pertain
to f\/lipping simultaneously two (disjoint) copies of
this edge in the double graph, f\/lip morphisms on pending edges (see Fig.~\ref{fi:mcg-pending}) change the topological
structure of the double graph (the genus and
the number of holes of the doubled Riemann surface may then change; only the total Poisson dimension must remain invariant).

As in the case of dual Teichm\"uller spaces in \cite{Fock2}, the coordinates $Z_\alpha$ describe the
linear subspace of the Teichm\"uller space ${\mathcal T}_{\hat g,\hat s}^H$, which is again a
$2^{\hat s}$-fold covering of the Teichm\"uller space ramif\/ied over
punctured surfaces (with old and new boundary components). This subspace comprises surfaces admitting the
involution interchanging two halves of the doubled Riemann surface.

\begin{figure}[tb]
\centerline{\includegraphics{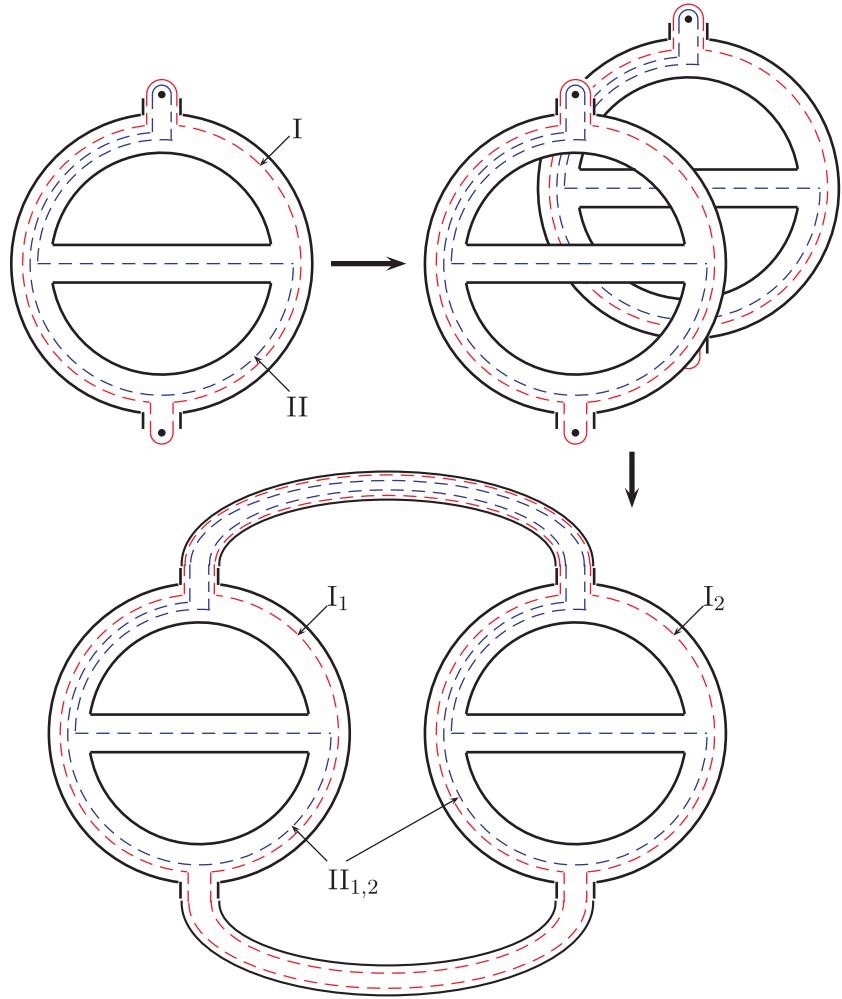}}
\caption{A general procedure of doubling a multicurve. We f\/irst produce two copies of the
original pattern, then sew together pairwise (with orientation preserving and with parallel transport of
the multicurve threads) copies of dead ends forgetting about dot-vertices. In the doubled surface,
we then obtain two copies (I$_1$ and I$_2$) of those original geodesic lines (I) that experienced even number of
inversions in the original pattern (although these copies are not parallel, their lengths coincide in
the new pattern), whereas each original geodesic (II) that experienced odd number of inversions, generates
one new geodesic (II$_{1,2}$) of {\em doubled} length in the doubled surface.}
\label{fi:multicurve-doubling}
\end{figure}

It seems that no closed geodesic line corresponds to an initial closed geodesic with
{\em odd} number of inversions (the corresponding line on the doubled Riemann surface is nonclosed).
To understand how to overcome this trouble and to formulate the general procedure of constructing
the ``doubled" multicurve, let us consider the example in Fig.~\ref{fi:multicurve-doubling}. There, we have
the multicurve with two components, geodesic curves I and II, on the windowed Riemann surface. Curve~I experiences
two inversions while curve~II has one inversion. To produce a multicurve on the doubled surface, we f\/irst
merely double the pattern (f\/irst stage in Fig.~\ref{fi:multicurve-doubling}) and then connect each pending edge of
the original graph with its copy on the clone of this graph in a natural way, that is, preserving the orientation and
without introducing intersections of threads of the multicurve (second stage in Fig.~\ref{fi:multicurve-doubling}).

We note, f\/irst, that the pattern in the doubled surface thus obtained is always a multicurve (since we did not introduce
any new intersection). We now turn to its content.

It is easy to see that a geodesic line~I with the geodesic function $G_{\mathrm I}$ and the length $\ell_{\mathrm I}$
on the original surface that has {\em even} number of inversions produces {\em two} disjoint nonself\/intersecting (albeit not parallel)
geodesic lines I$_1$ and I$_2$ on the doubled Riemann surface, and these two new lines has the respective geodesic
functions $G_{\mathrm I_1}$ and $G_{\mathrm I_2}$ and lengths $\ell_{\mathrm I_1}$ and $\ell_{\mathrm I_2}$ such that
\begin{gather}
\label{GI}
G_{\mathrm I_1}= G_{\mathrm I_2}= G_{\mathrm I},\qquad \ell_{\mathrm I_1}=\ell_{\mathrm I_2}=\ell_{\mathrm I}.
\end{gather}

On the contrary, for a geodesic line II with the geodesic function $G_{\mathrm {II}}$ and the length $\ell_{\mathrm {II}}$
that has an {\em odd} number of inversions, this geodesics produces {\em one} nonself\/intersecting geodesic line~II$_{1,2}$ on
the doubled Riemann surface with the geodesic function $G_{\mathrm {II}_{1,2}}$ and length $\ell_{\mathrm {II}_{1,2}}$,
and it is easy to see that they satisfy the relations
\begin{gather}
\label{GII}
G_{\mathrm {II}_{1,2}}= \left(G_{\mathrm {II}}\right)^2-2,\qquad \ell_{\mathrm {II}_{1,2}}=2\ell_{\mathrm {II}},
\end{gather}
that is, we have then a single geodesic of doubled length.

From (\ref{GI}) and (\ref{GII}) it follows that we must take as the characteristic of a multicurve, or
a rational lamination, the sum of lengths of the constituent geodesics.
We must be able to construct this sum from the geodesic function or from a quantum geodesic function ensuring the
positiveness property: the length must be nonnegative function in the classical case or positive-def\/inite operator
in the quantum case. This is ensured by the following construction proposed in~\cite{ChP}.

\begin{definition}\label{6-1}
The {\em proper length\/} ${\rm p.l.}(\gamma)\equiv \ell_\gamma$ of a closed curve
$\gamma$ in the classical or quantum case is
constructed from the quantum geodesic operator $G^\hbar_\gamma$ as
\begin{gather*}
{\rm p.l.}(\gamma)\equiv \ell_{\gamma}=2\lim_{n\to\infty}\frac1n
\log 2T_n(G^\hbar_\gamma/2),
\end{gather*}
where we take the principal branch of the logarithm and $T_n$ are Chebyshev's
polynomials (cf.~(\ref{cheb})). Since $T_n(\cosh {t\over 2})=\cosh {{nt}\over 2}$,
it follows that ${\rm p.l.}(\gamma )$ is the
hyperbolic length of $\gamma$ in the Poincar\'e metric in the classical case.

In the operatorial case, we can determine ${\rm p.l.}(\gamma )$ explicitly in terms of the
spectral expansion of the operator $G^\hbar_\gamma$, which is known exactly~\cite{Kashaev3}.
Namely, the basis of eigenfunctions of $G^\hbar_\gamma$ is
labeled by the positive number $S$ whereas the eigenvalue corresponding to the eigenfunction~$\left|\alpha _S\right\rangle$ has the form $\e ^{S/2}+\e ^{-S/2}$ and the both limits $S\to+\infty$ and $S\to 0$
are singular in the functional sense. Eigenfunctions $\left|\alpha _S\right\rangle$ constitute
an orthogonal and complete basis, so we can def\/ine the proper length
operator to be the one with the same eigenfunctions $\alpha _S$ and with positive
eigenvalues~$S/2$. The operator ${\rm p.l.}(\gamma )$ is then a
well-def\/ined operator on any compactum in the function space~$H$.

The {\em proper length\/} of a QMC or GMC $\hat{\mathrm C}$,
which we denote as ${\rm p.l.}(\hat C)\equiv \ell_{\hat C}$, is the sum of the proper lengths
of the constituent geodesic length operators
(or the sum of geodesic lengths calculated in the Poincar\'e metric
in the classical case) weighted by the number of appearances in the multiset.
\end{definition}

\begin{lemma}\label{lem-multicurve}
Given a general even-based multicurve ${\mathcal C}_e$, its {\em proper length}
\begin{gather}
\label{GMC-length}
\ell_{{\mathcal C}_e}\equiv\sum_{i=1}^n s_i\ell_{\gamma_i},
\end{gather}
satisfies the general relation: denoting by $\widehat {\mathcal C}_{1,2}$ the GMC generated on the doubled Riemann
surface, we have
\begin{gather*}
\ell_{\widehat {\mathcal C}_{1,2}}=2\,\ell_{{\mathcal C}_e}
\end{gather*}
irrespectively on the content of this multicurve: how many geodesic curves constitute the corresponding
multicurve, what is the number of inversions experienced by individual curves, etc.
\end{lemma}

We therefore see that the only object for which the Thurston theory can be elaborated are sums
(\ref{GMC-length}) of lengths of curves entering a multicurve. Note that we have
the convergency theorem~\cite{ChP} only for this characteristic.

But if we want to consider algebraic structures in a consistent way, it is better to do in the
original pattern (as was done above)
because geodesic algebras on a doubled surface do not satisfy the standard relations
due to the additional involutional symmetry.

\subsection{Elements of Thurston theory of bordered surfaces}\label{s:Thurston}

We now extend the notion of train tracks and foliation-shear coordinates to the case of bordered surfaces.
(In this section, we use the terminology of Thurston's school.) For a comprehensive research into
the subject, see \cite{HarPen}; an introductory description close to the one in the present paper
is contained in~\cite{ChP}.

\subsubsection{Decorated measured foliations and freeways}\label{ss:freeways}

\begin{figure}[tb]
\centerline{\includegraphics{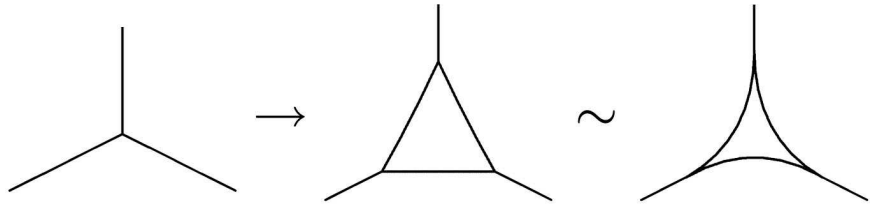}}
\caption{Freeway from fatgraph.}
\label{fi:freeway}
\end{figure}

We recall material from \cite{PP2}. If
$\Gamma_{g,\delta}\subseteq \Sigma_{g,\delta}$ is a cubic fatgraph spine of
$\Sigma_{g,\delta}$, then we may blow-up each three-valent vertex of
$\Gamma_{g,\delta}$ into a little trigon as illustrated in Fig.~\ref{fi:freeway}.  The resul\-ting object
$\tau=\tau_{\Gamma_{g,\delta}}$ has both a natural branched one-submanifold structure and a fattening, and
furthermore, components of $\Sigma_{g,\delta}\backslash\tau$ are either little trigons, or once-holed (once-punctured)
nullgons, or bigons with one marked point on one of the edges (we therefore forbid monogons).
Thus, $\tau$ is not a train track, but it is almost a train track, and is
called the {\it freeway} associated to $\Gamma_{g,\delta}$.  Notice that each edge of
$\Gamma_{g,\delta}$ gives rise
to a corresponding large branch of $\tau$, and each vertex gives rise to three small
branches.  It is easy to see that every {\em measured lamination} in
$\Sigma_{g,\delta}$ is carried by the freeway $\tau$.

The frontier of a once-holed
nullgon component of $\Sigma_{g,\delta}\backslash\tau$ is a puncture-parallel curve called a {\it collar
curve} of $\Sigma_{g,\delta}$ (curve~I in Fig.~\ref{fi:center}).
We also introduce collar curves that are unions of bigon edges constituting
a path homeomorphic to a boundary component (curve~II in Fig.~\ref{fi:center}). We have exactly
$s$ nonhomeomorphic types of collar curves.
A small branch is contained in exactly one collar curve, while a
large branch may be contained in either one or two collar curves (and each pending branch is contained
in exactly one collar curve).

A {\it measure} on a freeway $\tau$ is a function $\mu \in {\Bbb R}^{\# B(\tau )}$
($B(\tau )$ is the set of all branches, long and short)
satisfying the natural {\em switch conditions} provided that for each switch
$v$ of $\tau$, we have
\[
\sum_{{\rm outgoing}\atop~{\rm half-branches}~b}\mu (b)=
\sum_{{\rm incoming}\atop{\rm half-branches}~b}\mu (b).
\]

Such a function satisfying the switch conditions is commonly called a {\it
{\rm(}transverse\/{\rm)} measure} on $\tau$, and $\tau$ itself is said to be
{\it recurrent} train track if it supports a {\it positive} measure $\mu$
with $\mu (b)>0$ for each branch of $\tau$.

For freeways, however, the measure is not
necessarily nonnegative (as it is for train tracks).  Let $U(\tau )$ denote the
vector space of all measures on $\tau$.  Notice that $\mu\in U(\tau )$ is uniquely
determined by its values on the short branches (the switch conditions are
equivalent to the ``coupling equations''
\[
\mu (a_1)+\mu (b_1)=\mu (e)=\mu (a_2)+\mu (b_2),
\]
for any long branch $e$ whose closure contains the switches $v_1\neq v_2$ and
$a_i,b_i$ are the short branches incident on $v_i$ for $i=1,2$). And, simultaneously,
the values on the long branches alone also uniquely determine $\mu$. At the same time,
the {\em both} sets of transverse measures on short and long branches are subject to
constraints imposed by the positivity condition in the case of recurrent train tracks
(the triangle inequality in the case of measures on long branches), so the problem is to
propose new variables describing the transverse measures that are unconstrained.
We then identify $U(\tau )\approx{\Bbb R}^{LB(\tau )}$, where $LB(\tau )$
denotes the set of long branches of $\tau$.

Recall that, having nonnegative transverse measures on edges of a train track $\tau$,
we can canonically associate to $\tau$ the corresponding {\em measured foliation}
${\cal MF}_0(\Sigma_{g,\delta})$, which is determined up to the equivalence relation that
claims that if two foliations dif\/fer only by their
``collar weights'', i.e., the transverse measures of boundary-parallel arcs, then they are
equivalent, and in each equivalence class we f\/ind a unique measured foliation with all
collar weights equal to zero.

A nonnegative measure on $\tau$ then canonically
determines a point of ${\cal MF}_0(F)$ together with a~nonnegative ``collar
weight'', and we can def\/ine the space
$\widetilde{\cal MF}_0(F)={\cal MF}_0(F)\times{\Bbb R}^s$ of {\it decorated measured
foliations}.

We then have the theorem.

\begin{theorem}[\cite{PP2}]
The space $U(\tau)\approx {\Bbb R}^{LB(\tau )}$ gives global coordinates on
$\widetilde{\cal MF}_0$, and there is a~canonical fiber bundle
$\Pi :\widetilde{\cal MF}_0\to{\cal MF}_0$, where the fiber over a point is the set
${\Bbb R}^s$ of all collar weights on $F$.
\end{theorem}

\begin{remark}
{\rm
The natural action of mapping class group
$MC_g^s$ is by bundle isomorphisms of $\Pi$.  Furthermore,
$\Pi$ admits a natural $MC_g^s$-invariant section
$\sigma :{\cal MF}_0\to\widetilde{\cal MF}_0$
which is determined by the condition  of identically
vanishing collar weights.  The restriction of $\sigma$ to
${\cal MF}_0\subseteq\widetilde{\cal MF}_0$ gives a piecewise-linear embedding of the
piecewise-linear manifold ${\cal MF}_0$ into the linear manifold (vector space)
$\widetilde{\cal MF}_0\approx U(\tau)\approx {\Bbb R}^{LB(\tau )}$.
}
\end{remark}

\subsubsection{Shear coordinates for measured foliations}\label{scmf}

We now give an equivalent parametrization of measured foliations in terms of
``Thurston's shear coordinates'' that are closest analogues of
Thurston's shear coordinates $Z_\alpha$ on ${\cal T}^H(\Sigma)$.

We assign a corresponding signed quantity (positive for right, negative for left) as follows.
Given a measure $\mu$ on the long branches of the freeway $\tau$ associated to the fatgraph spine
$\Gamma\subseteq \Sigma$, def\/ine the {\it (Thurston's foliation-)shear coordinate} of the edge indexed by
$\alpha$ to be
\begin{gather}
\label{shear}
\zeta_\alpha ={1\over 2}(\mu (A)-\mu (B)+\mu (C)-\mu (D)),
\end{gather}
in the notation of Fig.~\ref{fi:flip} for nearby branches.

If $\alpha$th edge is a pending edge, we use also formula (\ref{shear}) but with, say, $\mu(C)=\mu(D)=0$
because the edge $\alpha$ has then just one switch $v_1$ at its upper end.

From the very def\/inition, $\zeta _\alpha$ is independent of collar weights.

The shear coordinates $\zeta _\alpha$ are not independent, but the only restrictions
imposed are
\begin{gather*}
\sum_{\alpha\in I}\zeta _\alpha=0
\end{gather*}
for the sums over edges $\alpha\in I$ surrounding any given boundary component;
we call these conditions the {\it face conditions} for shear coordinates.
Thus, the space of foliation-shear coordinates is of dimension
${LB(\tau)-s}$ and coincides, in particular, with the dimension of the nondegenerate Poisson leaf.
For any assignment of shear coordinates, there is a well-def\/ined point of ${\cal MF}_0$
realizing them, thereby
establishing a homeomorphism between ${\cal MF}_0$ and this sub-vector space ${\Bbb R}^{LB(\tau
)-s}\subseteq {\Bbb R}^{LB(\tau )}$ of shear coordinates on the long branches of $\tau$.

\begin{remark}\label{rem-fol-shear}
Note that for arbitrary integer measures $\mu(Z_\alpha)$ on long edges, the foliation-shear
coordinates $\zeta_\alpha$ on pending (long) edges can be half-integer! They are integers
only if we claim the evenness condition. In the picture of the double, these coordinates are
halves of the corresponding would-be foliation-shear coordinates $\zeta_\alpha$ on $\Sigma_{\hat g, \hat s}$ (where the
foliation is induced from the one in $\Sigma_{g,\delta}$ in a natural way described in
Fig.~\ref{fi:multicurve-doubling}). That we take halves of them follows from
the doubling of the corresponding shear coordinate $Z_\alpha$.

For example, the value of the foliation-shear coordinate on the pending edge for the pattern in
Fig.~\ref{fi:bounding} is $+2$; it is half-integer only if we have originally nonclosed curves.
\end{remark}

We now describe the action of the mapping class group (the Whitehead moves) on foliation-shear coordinates,
that is, derive the analogue of formula (\ref{abc}) for measured foliations,
which is an elementary calculation using the formulas for splitting.

\begin{figure}[tb]
\centerline{\includegraphics[scale=0.95]{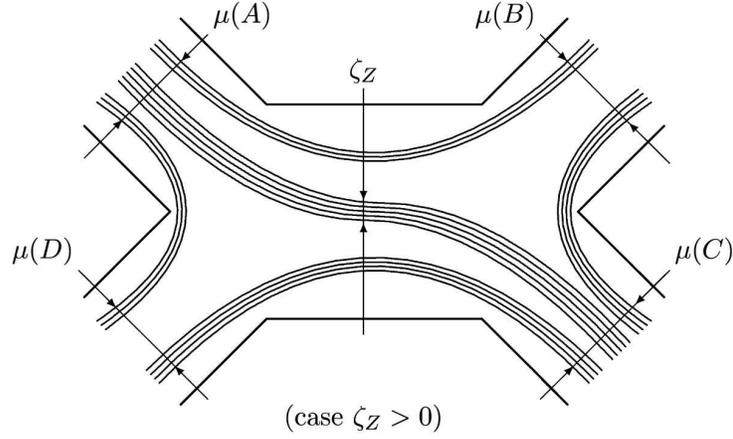}}
\caption{Foliation-shear coordinates.}
\label{fi:shear}
\end{figure}

\begin{figure}[t]
\centerline{\includegraphics{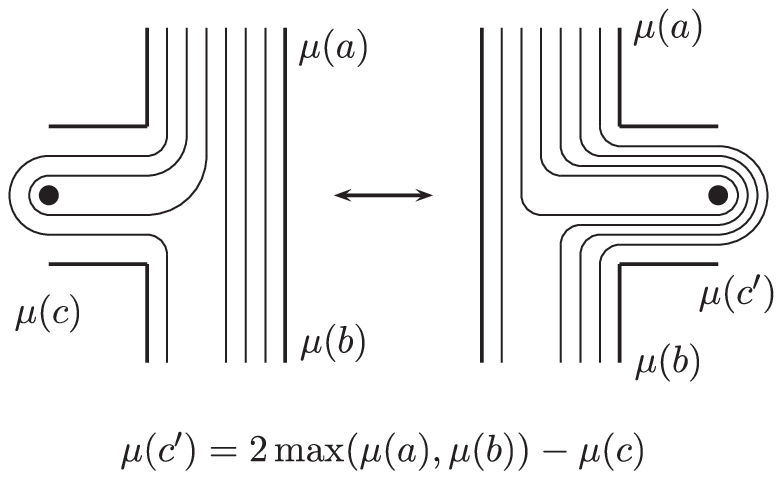}}
\caption{Mapping class group transformation of f\/lipping the pending edges acting on
transverse measures.}\vspace{-2mm}
\label{fi:fol-pend}
\end{figure}

\begin{lemma}\label{lem4-1}\looseness=-1
Under the Whitehead move in Fig.~{\rm \ref{fi:flip}}, the corresponding foliation-shear coordina\-tes
of the edges $A$, $B$, $C$, $D$, and $Z$ situated as in that figure are transformed
according to formula~\eqref{abc}
\begin{gather}
M_Z: (\zeta _A,\zeta _B,\zeta _C,\zeta _D,\zeta _Z)
\nonumber
\\
\qquad{}\mapsto(\zeta _A+\phi_H(\zeta _Z),
\zeta _B-\phi_H(-\zeta _Z), \zeta _C+\phi_H(\zeta _Z), \zeta _D-\phi_H(-\zeta _Z),-\zeta _Z)
\label{shear-flip}
\end{gather}
with
\begin{gather}
\label{h-trans}
\phi_H(\zeta_Z)=(\zeta_Z+|\zeta_Z|)/2,
\end{gather}
i.e., $\phi_H(x)=x$, for $x>0$, and zero otherwise.  All other
shear coordinates on the graph remain unchanged. Formula \eqref{shear-flip} holds
irrespectively on whether some of edges $A$, $B$, $C$, and $D$ are pending edges.

When flipping a pending edge as shown in Fig.~{\rm \ref{fi:mcg-pending}}, the corresponding
foliation-shear coordinates of the edges $Y_1$, $Y_2$, and $Z$ undergo the transformation
\begin{gather*}
M_Z: (\zeta _{Y_1},\zeta _{Y_2},\zeta _Z)\mapsto
(\zeta _{Y_1}-\phi_H(-2\zeta _Z), \zeta _{Y_2}+\phi_H(2\zeta _Z), -\zeta _Z)
\end{gather*}
with $\phi_H(x)$ from \eqref{h-trans}. The corresponding transformations for
the transverse measures are depicted in Fig.~{\rm \ref{fi:fol-pend}}.
\end{lemma}

\begin{remark}\label{rem4-1}
Comparing expressions for
the classical function $\phi(x)=\log(1+\e^x)$ and \eqref{h-trans}, one f\/inds that the
latter is a {\em projective limit\/} of the former:
\begin{gather*}
\phi_H(x)=\lim_{\lambda\to+\infty}\frac1\lambda \phi(\lambda x)
=\lim_{\lambda\to+\infty}\frac1\lambda \phi^\hbar(\lambda x),
\end{gather*}
that is, {all three} transformations coincide {asymptotically}
in the domain of large absolute values {\rm(}or large eigenvalues for the
corresponding operators\/{\rm)}
of Teichm\"uller space coordinates $\{Z_\alpha\}$.
This property was crucial for constructing quantum Thurston theory of surfaces
with holes in~\cite{ChP}. Presumably, the corresponding construction can be
transferred almost literally to the case of bordered surfaces. We do not however
consider it in this paper postponing the description of quantum Thurston theory
of bordered surfaces for future publications.
\end{remark}

\section{Conclusion}

We have proposed a (novel) graph description of moduli spaces of Riemann surfaces with windows
(marked points on the boundary) and demonstrated that it provides all the necessary ingredients, including the
construction of the double Riemann surface with holes (and without windows).

The approach in this paper might be close to the one in~\cite{FST} where mutations pushing marked points
from one boundary component to another have been considered as well. There, however, the corresponding points
can be transferred only in pairs, which apparently dif\/fers from the procedure in this paper. Establishing a
correspondence between these two approaches def\/initely deserves further studies.

Worth mentioning is the open/closed string diagrammatics proposed recently in \cite{KP}. We hope that
introducing new variables (on pending edges) may allow attaining a comprehensive quantitative description
of string theories and, in its possible quantum version, may provide a bridge to a string f\/ield theory description
of the related objects.

Also note the algebras $A_n$ and $D_n$, which play an important role in classif\/ication of f\/inite-type cluster
algebras \cite{FST}. Eventually, it seems interesting to generalize the construction in this paper to higher
Teichm\"uller spaces introduced by Fock and Goncharov \cite{FG}.

\subsection*{Acknowledgments}

The author is indebted to V.V.~Fock and R.C.~Penner for the fruitful discussion on the Oberwolfach Conference
on Teichm\"uller spaces, which initiated this work.

This work has been partially f\/inancially supported by the
RFBR Grant No.~05-01-00498, by the Grant for Support of the Leading Scientif\/ic
Schools 2052.2003.1, by the Program Mathematical Methods of
Nonlinear Dynamics, by the ANS Grant ``G\'eom\'etrie et Int\'egrabilit\'e en Physique Math\'ematique''
(contract number ANR-05-BLAN-0029-01), and
by the European Community through the FP6
Marie Curie RTN {\em ENIGMA} (Contract number MRTN-CT-2004-5652).

\pdfbookmark[1]{References}{ref}
\LastPageEnding

\end{document}